\documentclass[onecolumn,10pt]{article}
\pdfoutput=1

\usepackage{PRIMEarxiv}

\usepackage[utf8]{inputenc} 
\usepackage[T1]{fontenc}
\usepackage{lmodern}
\usepackage{url}            
\usepackage{booktabs}       
\usepackage{amsfonts}       
\usepackage{nicefrac}       
\usepackage{microtype}      
\usepackage{lipsum}
\usepackage{fancyhdr}       
\usepackage{graphicx}       

\usepackage{mathtools,amssymb}

\usepackage{setspace}
\usepackage{algorithmic}
\usepackage{algorithm}

\usepackage{amsmath}
\usepackage{amssymb}
\usepackage{mathtools}
\usepackage{amsthm}
\usepackage{amsfonts}

\newcommand{\argmin}{\operatornamewithlimits{argmin}}

\newcommand{\minimize}{\operatornamewithlimits{minimize}}

\newcommand{\subt}{\operatornamewithlimits{\text{subject to}}}

\usepackage{subcaption}

\usepackage{gensymb}

\usepackage[capitalize,noabbrev]{cleveref}

\graphicspath{{media/}}

\usepackage{titling}
\thanksmarkseries{alph}

\title{Maximum Likelihood Localization of a Network of Moving Agents from Ranges, Bearings and Velocity measurements}

\author{
	Filipa Valdeira \thanks{Computer Science Department, NOVA LINCS, NOVA School of Science and Technology, Universidade NOVA de Lisboa, Portugal} \thanks{Center for Mathematics and Applications, NOVA Math, NOVA School of Science and Technology, Universidade NOVA de Lisboa, Portugal } \thanks{Corresponding author (email:\texttt{f.valdeira@fct.unl.pt})}
	\and  
	Cl\'audia Soares \thanksmark{1} 
	\and 
	Jo\~{a}o Gomes\thanks{Institute for Systems and Robotics (ISR/IST), LARSyS, Intituto Superior Técnico, Universidade de Lisboa, Portugal} 
}

\date{}

\graphicspath{{./Figures/}}

\begin{document}
	\maketitle
	\begin{abstract}
	Localization is a fundamental enabler technology for many applications, like vehicular networks, IoT, and even medicine. While Global Navigation Satellite Systems solutions offer great performance, they are unavailable in scenarios like indoor or underwater environments, and, for large networks, the instrumentation cost is prohibitive. We develop a localization algorithm from ranges and bearings, suitable for generic mobile networks.  Our algorithm is built on a tight convex relaxation of the Maximum Likelihood position estimator. To serve positioning to mobile agents, a horizon-based version is developed accounting for velocity measurements at each agent. To solve the convex problem, a distributed gradient-based method is provided. This constitutes an advantage over centralized approaches, which usually exhibit high latency for large networks and present a single point of failure. Additionally, the algorithm estimates all required parameters and effectively becomes parameter-free. Our solution to the dynamic network localization problem is theoretically well-founded and still easy to understand. We obtain a parameter-free, outlier-robust and trajectory-agnostic algorithm, with nearly constant positioning error regardless of the trajectories of agents and anchors, achieving better or comparable performance to state-of-the-art methods, as our simulations show. Furthermore, the method is distributed, convex and does not require any particular anchor configuration.
	\end{abstract}
	\keywords{Maximum Likelihood estimation, Dynamic network localization, Hybrid measurements, Convex optimization, Distributed optimization}

	\section{Introduction}
	
While multi agent systems often rely on Global Navigation Satellite Systems (GNSS) to accomplish their tasks, some environmental constraints prevent their use. For example, for Autonomous Underwater Vehicles (AUVs) \cite{article:Miller10,article:Yoo14}, where the high conductivity in water decreases radio wave propagation, hence preventing the use of satellite positioning. This paper addresses the network localization problem for GNSS-denied environments from ranges, bearings, and velocity noisy measurements, {with emphasis on AUV applications}. Advancements in technology have made communication between these vehicles possible, so that agents may share information to obtain better localization accuracy. This is called a \textit{cooperative} approach, which is also a trend in other GNSS-deprived settings, such as indoor environments \cite{article:UAV_coop,article:coop_indoors}. Besides, as
networks grow in size, \textit{distributed} approaches (as opposed to centralized) are often desirable. We define distributed operation of a network of sensing and computing nodes as a mode where each agent performs the same computations and only communicates with one-hop neighbors on the communication graph. As the distributed system does not centralize on a single processor all available information, it avoids the communication strain near the central node and on the bottleneck hubs (observed in centralized approaches), thus imposing a small latency. Importantly, the distributed paradigm does not present a single point of failure, traditionally found on federated star topologies. Therefore, we target the case of generic agent networks, able to cooperate with each other to accomplish self-localization tasks in a distributed manner, while performing their own trajectories, i.e., in a dynamical setting.

\paragraph*{Problem formulation} A common approach for position estimation in dynamic settings is the use of Bayesian Filters, such as the Extended Kalman Filter (EKF), which is very often employed for localization tasks
\cite{article:Yoo14,article:Miller10,article:Pinheiro17,article:Dong12,article:EKF_2023}. In spite of its popularity, EKF introduces linearization errors (as nonlinear systems must be linearized) and requires correct initialization to avoid convergence problems, which may not be available. In an attempt to overcome these challenges, we rely instead on the Maximum Likelihood Estimator (MLE) 	which is, for large sets of data and high signal to noise ratio, asymptotically unbiased and asymptotically efficient \cite{book:Kay97}. Therefore, the advantage over competing formulations is that if the error does follow the considered distribution, MLE will be the optimal estimator.

\paragraph*{Nonconvexity and relaxations}  
Regardless of the formulation, the localization problem typically leads to a \textit{nonconvex} cost function. While it is possible to preserve the nonconvex function \cite{article:calafiore10, article:distributed_nonconvex_costa2006,article:Erseghe2015}, this requires a correct initialization of the method, which can be difficult to obtain. {In order to} avoid this limitation, we aim at a convex formulation, which under the MLE formulation calls for a relaxation of the original problem. Popular relaxation strategies for the nonconvex MLE include Semidefinte Programming (SDP) \cite{article:Biswas06,YANG2022108504} and Second-Order Cone Programming (SOCP) \cite{article:SOCP_Tseng_2006}. {Both approaches have strong requirements on the network topologies, in particular regarding anchor placement. For example, SOCP relaxations require the nodes to remain in the convex hull of the anchors, which can be a limitation when the area of actuation in unknown a priori. Other proposals include a Sum of Squares relaxation after problem reformulation \mbox{\cite{article:SOS_2006}}, which is a computationally expensive alternative. Our proposed inclusion of angle measurements in the position estimation proves helpful in removing these constraints on network configuration and consequentially allows more flexibility in trajectories, while remaining suitable for large networks.}

\paragraph*{Measurement availability} While the previous methods consider only pairwise range measurements between agents \cite{article:Biswas06,article:distributed_SOCP,article:Shi10,SOARES2021}, pairwise bearings may also be available, leading to bearing-only \cite{article:Wang23} or hybrid methods \cite{article:Crouse13,article:convex_optimization_uav,article:FLORIS,article:Biswas5,article:SDP_benchmark}. Here, we consider the hybrid setting, which has been shown to be beneficial when compared to the single measurement case \cite{article:Hybrid_Eren2011}. We consider that all nodes have access to range measurements with respect to their neighbours and a subset of them are able to additionally retrieve bearings {from a different device}. This differs, for instance, from \cite{article:mixed_Lin} where each node has access to a single type of measurement, which can be of a different nature. We note that one important line of work takes into account Radio Signal Strength (RSS) and direction measurements \cite{article:HybridRSS_Ding2021,article:hybrid_RSSAoA_Tomic2019}, which calls for considerably different range models, than those assumed here.

\paragraph*{ML Hybrid localization} Similarly to range-only, we find ML \cite{article:Crouse13,article:SDP_benchmark} and non-ML \cite{article:convex_optimization_uav,article:FLORIS,article:Biswas5} derived cost functions for the hybrid setting. In~\cite{article:Crouse13} the authors compare different distributions of angle measurements, where ranges are often available. The formulation is closely related to ours, although derived for a 2D case and for single target location. The nonconvexity is not relaxed and a Least Squares approach is used for initialization, after which the optimization is carried out with gradient descent. Consequently, this is subject to local minima and under higher noise levels some points are poorly estimated.{ In~\mbox{\cite{article:SDP_benchmark}} and \mbox{\cite{article:Soares20}}, the authors follow a ML formulation where directional data is modeled with the von Mises-Fisher distribution. Our proposed method uses the relaxation in \mbox{\cite{article:Soares20}} that is found to outperform the SDP relaxation \mbox{\cite{article:SDP_benchmark}} in positioning accuracy and is more amenable to larger networks \mbox{\cite{article:Soares20}}.}

\paragraph*{Mobile networks} While methods devised for static networks can be applied for position estimation in mobile networks, it is beneficial to look at different time instants to retrieve additional information. An example is the bearing-only approach proposed in \cite{article:Wang23} where velocity is estimated alongside position. Other settings assume that relative velocity measurements between nodes are available \cite{article:Kumar16,article:RangeVelRelative} and extend existing range-only approaches. Instead, we assume that each node has access to its own velocity measurements (in addition to range and bearing measurements). Thus, we extend the hybrid ML formulation to account for individual velocity measurements, adding a time-window to account for different time instants. We show that this formulation is still amenable to the relaxation proposed in \cite{article:Soares20}, thus overcoming the limitations in alternative relaxations to the hybrid ML problem. However, we emphasize that we assume no dynamical model for the movement, unlike the standard EKF approach. 

\paragraph*{Distributed setting} Given their centralized nature, most of the previous approaches are not suitable for networks with a large number of nodes, which are often found in current applications. Consequently, there has been a growing interest in distributed methods \cite{article:distributed_SOCP,article:Shi10,article:distributed_DILAND,article:Simonetto,article:Wang23, article:DistribCoop_Salari,article:KalmanRangeBear}. In particular, distributed versions of the SOCP relaxation \cite{article:distributed_SOCP} and of the edge-based SDP relaxation \cite{article:Shi10} have been proposed for the range-only setting, {but they still imply strong assumptions on network configuration}. We show that our dynamic formulation for the hybrid case upon convex relaxation can be solved in a distributed manner. Thus, we reach a convex method (not dependent on initialization) which benefits from hybrid pairwise measurements (improving positioning accuracy) and takes into account dynamic information (without specifying an underlying dynamical model), while being scalable to large networks of agents (due to its distributed implementation).

\paragraph*{Noise parameters.} We note that, while the convex relaxation of the MLE does not need a precise initialization, it does require knowledge of noise parameters beforehand. This aspect is often overlooked and it is assumed that these values are known or previously obtained by a different method \cite{article:SDP_benchmark,article:Biswas5,article:convex_optimization_uav}, which limits straightforward applicability to new scenarios. In contrast, our approach {includes the estimation of noise distribution parameters}, making it parameter-free and more versatile.

\subsection{Contributions} We extend the previous work through three main contributions:
\begin{itemize}
	\item We propose a horizon-based version of the hybrid MLE formulation, accounting for velocity measurements at each agent, which is better suited for mobile agents but still amenable to a tight convex relaxation (Section~\ref{sec:dynamic}). The final algorithm allows for a more robust estimation of position compared to existing approaches (Section~\ref{sec:simulation}).
	\item We obtain a distributed algorithm that makes the algorithm scalable for large networks and does not lead to a single point of failure (Section~\ref{sec:imple}).
	\item We derive a parameter-free version of the method with a minimal decrease in positioning accuracy (Section~\ref{sec:parameters}), that can be easily implemented without \textit{a priori} knowledge.
\end{itemize}

\section{Problem Formulation and Background}
\label{sec:form}
In light of previous considerations, we present a localization algorithm for a network of a variable and undefined number of moving autonomous agents. Anchors are assumed to be vehicles with access to an accurate positioning system, either by benefiting from GNSS signals or by remaining fixed at a known position. Besides, each vehicle is assumed to have a ranging device providing noisy distance measurements between itself and other vehicles within a certain range. Some vehicles are also equipped with a vector sensor producing bearing measurements. This angular information is assumed to be already in a common frame of reference, which could be obtained resorting to compass measurements, so that all angles are given with respect to the North direction, for instance. It is further considered that if vehicle ``A'' has measurements with respect to vehicle ``B'', then the reverse is also true. All vehicles are potentially mobile and can follow any trajectory, including anchors.

The main goal is, then, to determine the position of each autonomous agent on a network, during a certain time interval. For this purpose, the network is represented as an undirected connected graph $\mathcal{G}(t)=(\mathcal{V}(t),\mathcal{E}(t))$, where each different graph corresponds to a given time instant. The vertices, also called nodes, $\mathcal{V}=\{1,...,n\}$ correspond to the vehicles, while the edges $i \sim j \in \mathcal{E}^x(t) $ indicate the existence of a measurement of type $x$ between node $i$ and $j$, at time $t$. In particular, $\mathcal{E}^d(t) $ refers to distance measurements and $\mathcal{E}^u(t) $ to angle measurements between two nodes. The known vehicle positions required for localization (anchors) are indicated by the set $\mathcal{A}=\{1,...,m\}$. Furthermore, the subset $\mathcal{A}^x_i(t) \subset \mathcal{A}(t)$ contains the anchors relative to which node $i$ has an available measurement of type $x$ at time $t$. Specifically, $\mathcal{A}^d_i $ refers to distance measurements and $\mathcal{A}^u_i $ to angle measurements to anchors.

We formulate our localization problem in $\mathbb{R}^p$ space. Practical relevant cases are $p=2$ for a bidimensional problem and $p=3$ for a tridimensional one. The position of node $i$ at time $t$ is indicated as $x_i(t) \in \mathbb{R}^p$, while anchor positions are designated as $a_k(t) \in \mathbb{R}^p$. Each node collects measurements with respect to its neighboring nodes or nearby anchors within a certain distance. Noisy distance measurements at time $t$ are represented as $d_{ij}(t) = d_{ji}(t)$ if they exist between two nodes ($i$ and $j$) or as $r_{ik}(t)$ if they occur between a node $i$ and anchor $k$.  Noisy angle measurements between node $i$ and $j$ are represented as a unit-norm vector in the correspondent direction $u_{ij}(t)$. Similarly, noisy angle measurements between node $i$ and anchor $k$ are taken as the unit-norm vector $q_{ik}(t)$. 

The problem is, thus, to estimate the unknown vehicle positions $x(t) = \{x_i(t) : i \in \mathcal{V} \}$ given a dataset of measurements
\begin{equation*} \label{eq:data}
	\begin{split}
		\mathcal{D}(t) =& \underbrace{\{ d_{ij}(t) : i \sim j \in \mathcal{E}^d(t) \}}_{\text{node-node ranges}} \bigcup \underbrace{\{ r_{ik}(t)
			: i \in \mathcal{V}, k\in \mathcal{A}^d_i(t) \}}_{\text{node-anchor ranges}} \\ &\bigcup \underbrace{\{ u_{ij}(t) :
			i \sim j \in \mathcal{E}^u(t) \}}_{\text{node-node angles}} \bigcup \underbrace{\{ q_{ik}(t) : i \in
			\mathcal{V}, k\in \mathcal{A}^u_i(t) \}}_{\text{node-anchor angles}}.
	\end{split}
\end{equation*}

A common assumption is to model distance noise with a Gaussian distribution of zero mean and variances $\sigma_{ij}^2$ and
$\varsigma_{ik}^2$, for node-node and node-anchor edges, respectively. Bearing noise is better represented by a von
Mises-Fisher distribution, specifically developed for directional data, with mean direction zero and concentration parameter
$ \varkappa_{ij}$ and $\lambda_{ik}$, for node-node and node-anchor edges, respectively. {Since we assume that nodes have two distinct devices for bearing and range measurements, we assume that the noise is uncorrelated between them. Therefore,} with the additional assumption of independent and identically distributed noise, the maximum likelihood estimator is given by the solution of the optimization problem
\begin{equation}
	\label{eq:MLE_optimization_problem}
	\begin{split}
		\minimize_{x} &\quad f_{\textrm{dist}}(x(t)) + f_{\textrm{ang}}(x(t))\quad ,
	\end{split}
\end{equation} where the functions are given as
\begin{equation}
	\begin{split}
		f_{\textrm{dist}}(x(t)) &=\sum_{i\sim j \in \mathcal{E}^d(t)} \frac{1}{2\sigma_{ij}^2}\big(\|x_i(t)-x_j(t)\|-d_{ij}(t)\big)^2  + \\ & \sum_{i \in \mathcal{V}} \sum_{k\in \mathcal{A}^d_i(t)}\frac{1}{2\varsigma_{ik}^2}\big(\|x_i(t)-a_k(t)\|-r_{ik}(t)\big)^2
	\end{split}
\end{equation}
for range measurement terms, and for the angular terms
\begin{equation}
	\begin{split}
		f_{\textrm{ang}}(x(t)) &=  -\sum_{i \sim j \in \mathcal{E}^u(t)} \Bigg( \varkappa_{ij} u_{ij}^T(t) \frac{x_i(t)-x_j(t)}{\|x_i(t)-x_j(t)\|} \Bigg) - \\ &
		\sum_{i \in \mathcal{V}}\sum_{k\in \mathcal{A}^u_i(t)}\Bigg( \lambda_{ik} q_{ik}^T(t) \frac{x_i(t)-a_k(t)}{\|x_i(t)-a_k(t)\|}  \Bigg)
	\end{split}\quad .
\end{equation}
Solving~\eqref{eq:MLE_optimization_problem} is difficult because the problem is nonconvex and nondifferentiable, especially for a large-scale network of agents. Problem \eqref{eq:MLE_optimization_problem} is nonconvex in both the distance and angle terms. {Nonconvexity in distance terms arises from the negative argument of the square for position differences smaller than $d_{ij}$, i.e., \mbox{$\|x_i-x_j\|<d_{ij}$}, while the angle terms include a non-linear denominator \mbox{$\|x_i-x_j\|$}. For further insights on the relaxation employed we refer the reader to \mbox{\cite{article:Soares20}}}.

\subsection{Convex Relaxation}
\label{sec:backg}
{We follow the relaxation proposed in \hbox{\cite{article:Soares20}} to deal with the nonconvexity of \hbox{Problem~\eqref{eq:MLE_optimization_problem}}. We start by summarizing this approach,} before extending it with velocity measurements in Section~\ref{sec:dynamic}. It should be noted that, throughout this section, we drop the time dependency notation for clarity. First, a new variable $y_{ij}$ is introduced to obtain the following equivalent formulation \begin{equation}
	\label{original_func}
	(\| x_i- x_j\|-d_{ij})^2 = \inf_{\|y_{ij}\|=d_{ij}} \|x_i-x_j-y_{ij}\|^2.
\end{equation} The constraint is then relaxed as $\|y_{ij}\|\leq d_{ij}$, with the same reasoning applied to node-anchor edges.

However, the problem is still not convex due to the angular terms. Consequently, they are relaxed with their proxies from the distance term minimization. Therefore, the denominator $\|x_i-x_j\|$ is approximated by $d_{ij}$ because this is exactly the purpose of the distance terms in the cost function. With similar reasoning, $x_i-x_j$ is approximated by $y_{ij}$. Summing up, the nodes' angular terms become $\frac{x_i-x_j}{\|x_i-x_j\|} \approx \frac{y_{ij}}{d_{ij}}$ and respectively for anchors $\frac{x_i-a_k}{\|x_i-a_k\|} \approx \frac{w_{ik}}{r_{ik}}$ . The final expression for this convex approximation is then
\begin{equation}
	\begin{split}
		\minimize_{x,y,w}&\quad f_{\textrm{dist}}(x,y,w)+f_{\textrm{ang}}(y,w)
		\\ \subt & \quad	\|y_{ij}\| \leq d_{ij}, \|w_{ik}\| \leq r_{ik}
	\end{split}\quad ,
	\label{convex approximation}
\end{equation}
where
\begin{equation*}
	\begin{split}
		f_{\textrm{dist}}(x,y,w) = &\sum_{i\sim j \in \mathcal{E}^d} \frac{1}{2\sigma_{ij}^2}\|x_i-x_j-y_{ij}\|^2 + \sum_{i \in \mathcal{V}} \sum_{k\in \mathcal{A}^d_i}\frac{1}{2\varsigma_{ik}^2}\|x_i-a_k-w_{ik}\|^2
	\end{split}
\end{equation*}
and 
\begin{equation*}
	\begin{split}
		f_{\textrm{ang}}(y,w)& =  -\sum_{i \sim j \in \mathcal{E}^u} \Bigg( \varkappa_{ij} u_{ij}^T \frac{y_{ij}}{d_{ij}} \Bigg)-\sum_{i \in \mathcal{V}} \sum_{k\in \mathcal{A}^d_i}\Bigg( \lambda_{ik} q_{ik}^T \frac{w_{ik}}{r_{ik}}  \Bigg)
	\end{split}\quad .
\end{equation*} This is a tight relaxation with respect to the original problem that outperforms comparable SDP state-of-the-art relaxations, as shown in \cite{article:Soares20}.
	
	\section{Dynamic Formulation}
	\label{sec:dynamic}
	
	While the previous formulation is adequate for static networks, when vehicles perform a trajectory, it is beneficial to take into account different consecutive positions to obtain a better estimate. Therefore, we introduce a time window and consider that each vehicle $i$ has access to its own velocity $\beta_i(t)$ at time $t$.
	
	In order to model velocity measurements, we decompose them into a unit vector $v_i(t)$ and magnitude $V_i(t)$. Under this format, we can look at $V_i(t)\Delta T$, the distance travelled by one node during a time interval $\Delta T$, as the distance between two points $x_i(\tau)$ and $x_i(\tau-1)$, where time has been discretized and is now indexed by an integer $\tau$. Besides, the unit vector $v_i(t)$ can be thought of as the bearing between those same nodes, $x_i(\tau)$ and $x_i(\tau-1)$, usually denoted as heading. Essentially, this means that we have distance and bearing measurements just as in the previous formulation, except that they exist between the same node at different time instants. So, we can model them in a similar way, with $V_i(t)\Delta T$ following a Gaussian distribution with variance $\sigma_i^2$ and $v_i(t)$ following a von Mises-Fisher distribution with concentration parameter $\varkappa_{i}$.
	
	Therefore, we obtain a formulation with similar structure to Problem~\eqref{eq:MLE_optimization_problem}, which can be relaxed in the same manner, by introducing an additional variable $s$, with the same role as $y$. We define the concatenation of the optimization variables over nodes or edges at instant $\tau$ as $x( \tau)=\{x_{i}(\tau)\}_{i \in \mathcal{V}}$, $y( \tau)=\{y_{ij}(\tau)\}_{i\sim j}$, $s( \tau)=\{s_{i}(\tau)\}_{i \in \mathcal{V}}$, $w_{i}(\tau)=\{w_{ik}(\tau)\}_{k\in \mathcal{A}_i}$,  $w(\tau)=\{w_{i}(\tau)\}_{i \in \mathcal{V}}$. Similarly, the concatenation over a time window of size $T_0$ at time $t$ is given as $x =\{x(\tau)\}_{t-T_0\leq \tau \leq t} $,  $y = \{y(\tau)\}_{t-T_0\leq \tau \leq t}$ , $s =\{s(\tau)\}_{t-T_0\leq \tau \leq t}$ and $w =\{w(\tau)\}_{t-T_0\leq \tau \leq t} $. Note that the dependency on $t$ (the time instant) is omitted from the notation in all variables ($x,y,s,w$) for compactness. 
	The final problem is given as
	\begin{equation}
		\label{eq:no_matrix_vel}
		\begin{split}
			\minimize_{x,y,w,s}&\quad f_{\textrm{dist}}(x,y,w)+ f_{\textrm{ang}}(y,w)+f_{\textrm{vel}}(x,s)
			\\ \subt & \quad	\|y_{ij}(\tau)\| \leq d_{ij}(\tau),\quad \|w_{ik}(\tau)\|\leq r_{ik}(\tau), \quad \|s_i(\tau)\| \leq V_i(\tau) \Delta T,
		\end{split}
	\end{equation} 
	where	
	\begin{equation}
		\label{vel_terms}
		\begin{split}
			f_{\textrm{vel}}&(x,s) = \sum_{i \in \mathcal{V} }\sum_{\tau =t-T_0+1}^{t}\frac{1}{2\sigma_{i}^2}\big(\|x_i(\tau) -x_i(\tau-1)-s_{i}(\tau)\|\big)^2-\\& \sum_{ i \in \mathcal{V}}\sum_{\tau =t-T_0+1}^{t} \varkappa_i \Bigg(  v_{i}(\tau)^T \frac{s_i(\tau)}{V_i(\tau) \Delta T} \Bigg)\quad, 
		\end{split}
	\end{equation}	
	\begin{equation}
		\label{dist_term}
		\begin{split}
			f_{\textrm{dist}}&(x,y,w)=	\sum_{i\sim j \in \mathcal{E}^d}\sum_{\tau =t-T_0}^{t}  \frac{1}{2\sigma_{ij}^2}\big(\|x_i(\tau) -x_j(\tau)-y_{ij}(\tau)\|\big)^2
			+\\& \sum_{i \in \mathcal{V}} \sum_{k\in \mathcal{A}^d_i}\sum_{\tau =t-T_0}^{t}  \frac{1}{2\varsigma_{ik}^2}\big(\|x_i(\tau)-a_k(\tau)-w_{ik}(\tau)\|\big)^2
		\end{split}
	\end{equation}
	and
	\begin{equation}
		\label{angle_terms}
		\begin{split}
			f_{\textrm{ang}}&(y,w) = -\sum_{i \sim j \in \mathcal{E}^u}\sum_{\tau =t-T_0}^{t}  \Bigg( \varkappa_{ij} u_{ij}(\tau)^T \frac{y_{ij}(\tau)}{d_{ij}(\tau)} \Bigg) - \\& \sum_{i\in \mathcal{V}} \sum_{k\in \mathcal{A}^u_i}\sum_{\tau =t-T_0}^{t}\Bigg( \lambda_{ik} q_{ik}(\tau)^T \frac{w_{ik}(\tau)}{r_{ik}(\tau)}  \Bigg).
		\end{split}
	\end{equation}
	
	\subsection{Reformulation}\label{sec:Formulation_final}
	
	We shall now reformulate \eqref{eq:no_matrix_vel} in matrix notation. Let $\Sigma_N$ be the diagonal matrix of $\frac{1}{\sigma_{ij}}$. Then, the first term of \eqref{dist_term} can be written as 
	\begin{equation}
		\label{eq:distance term M}
		f_{\textrm{dist,nodes}}(x,y) =   \frac{1}{2}\|\Sigma_N Ax - \Sigma_N y \|^2,
	\end{equation}
	where matrix $A$ is the Kronecker product of the identity matrix of dimension $T_0$ with the result from the Kronecker product of $C$, the arc-node incidence matrix of the network, with identity matrix of dimension $p$, that is,
	\begin{equation*}
		A= I_{T_0} \otimes(C \otimes I_p) = (I_{T_0} \otimes C) \otimes I_p.
	\end{equation*}
	Essentially, this extends matrix $C$ along the dimension of the problem ($p$) and time window ($T_0$), assuming the set of edges remains constant over $T_0$. 
	
	In the same way, taking $\Sigma_V$ to be the diagonal matrix of $\frac{1}{\sigma_{i}}$, the first term of \eqref{vel_terms} may be written as 
	\begin{equation}
		\label{eq:dist_vel}
		f_{\textrm{dist,vel}}(x,s) =   \frac{1}{2}\|\Sigma_V Nx - \Sigma_V s \|^2
	\end{equation}
	where $N= C_{\textrm{vel}} \otimes I_p$. However, it should be clearly noted that, in this formulation, the set of edges has changed, hence the different designation $C_{\textrm{vel}}$. Each node at time $t$ has an edge with its position at $t-1$ and $t+1$, except for the first and last instants of the time window.
	The last distance term, concerning anchor-node measurements, may be reformulated as
	\begin{equation}
		\label{dist_anchor}
		f_{\textrm{dist,anchors}}(x,w) = \frac{1}{2}\|\Sigma_A Ex-\Sigma_A \alpha - \Sigma_A w\|^2
	\end{equation}where vectors for anchor positions are concatenated in the usual way as $\alpha_{i}(\tau)=\{a_{ik}(\tau)\}_{k\in \mathcal{A}_i}$,  $\alpha(\tau)=\{\alpha_{i}(\tau)\}_{i \in \mathcal{V}}$, $\alpha =\{\alpha(\tau)\}_{t-T_0\leq \tau \leq t} $. $E$ is a selector matrix, indicating which node has a measurement relative to each anchor and $\Sigma_A $ the diagonal matrix of $\frac{1}{\varsigma_{ik}}$.
	
	Finally, the angle measurements are grouped as $\tilde{u}_{ij}(\tau) = \frac{\varkappa_{ij}u_{ij}(\tau)}{d_{ij}(\tau)}$, $\tilde{q}_{ik}(\tau) = \frac{\lambda_{ik}q_{ik}(\tau)}{r_{ik}(\tau)}$ and $\tilde{v}_{i}(\tau) = \frac{\varkappa_{i}v_{i}(\tau)}{V_{i}(\tau)\Delta T}$. Following the same reasoning as before, $u(\tau)=\{\tilde{u}_{ij}(\tau)\}_{i\sim j}$, $v(\tau)=\{\tilde{v}_{i}(\tau)\}_{i \in \mathcal{V}}$ and $q_{i}(\tau)=(\tilde{q}_{ik}(\tau))_{k\in \mathcal{A}_i}$,  $q(\tau)=\{q_{i}(\tau)\}_{i \in \mathcal{V}}$. Then, concatenating along the time window, $u= \{u(\tau)\}_{t-T_0\leq \tau \leq t} $ , $v =\{v(\tau)\}_{t-T_0\leq \tau \leq t} $ and $q =\{q(\tau)\}_{t-T_0\leq \tau \leq t} $. It follows that
	\begin{equation}
		\label{eq:ang_node}
		f_{\textrm{ang}}(y,w,s) = - u^T y- q^Tw - v^Ts.
	\end{equation}
	Considering the definitions in \eqref{eq:distance term M}, \eqref{eq:dist_vel}, \eqref{dist_anchor} and \eqref{eq:ang_node}, Problem \eqref{eq:no_matrix_vel} is reformulated as
	\begin{equation}
		\label{eq:ref_quad}
		\begin{split}
			\minimize_{x,y,w,s}& \quad f_{\textrm{dist,nodes}}(x,y) + f_{\textrm{dist,anchors}}(x,w)+f_{\textrm{dist,vel}}(x,s) + f_{\textrm{ang}}(y,w,s)
			\\ \subt & \quad	\|y_{ij}(\tau)\| \leq d_{ij}(\tau),\quad \|w_{ik}(\tau)\|\leq r_{ik}(\tau),\quad \|s_i(\tau)\| \leq V_i(\tau) \Delta T.
		\end{split}
	\end{equation}	
	Introducing a variable $z=(x,y,w,s)$ and defining $\mathcal{Z}=\{ z : \|y_{ij}(t)\| \leq d_{ij}(t), i\sim j \in \mathcal{V}; \|w_{ik}(t)\| \leq r_{ik}(t), i\in \mathcal{V}, k\in \mathcal{A}_i; \|s_i(t)\| \leq V_i(t) \Delta T , i\in \mathcal{V}\}$, Problem~\eqref{eq:ref_quad} can be written in a quadratic form as
	\begin{equation}
		\label{Final_Formulation_Z}
		\begin{aligned}
			& \minimize_{z}\quad \frac{1}{2}z^TMz-b^Tz\\
			& \subt \quad z \in \mathcal{Z}	
		\end{aligned}		
	\end{equation}
	where $M=M_1+M_2+M_3$, with 
	\begin{equation*}
		M_1=	\begin{bmatrix}A^T\Sigma_N  \\- \Sigma_N  \\ 0\\0\end{bmatrix}
		\begin{bmatrix}\Sigma_N A & - \Sigma_N  & 0& 0\end{bmatrix},\quad		M_2=\begin{bmatrix}E^T\Sigma_A   \\ 0 \\- \Sigma_A \\0\end{bmatrix}
		\begin{bmatrix}\Sigma_A E  & 0 &- \Sigma_A &0\end{bmatrix},	
	\end{equation*}
	\begin{equation*}
		M_3=	\begin{bmatrix}N^T\Sigma_V  \\0 \\ 0\\- \Sigma_V \end{bmatrix}
		\begin{bmatrix}\Sigma_VN & 0 & 0& - \Sigma_V\end{bmatrix},	
	\end{equation*}
	and $b = b_1+b_2$, with 
	\begin{equation*}
		b_1 = \begin{bmatrix}E^T\Sigma_A   \\ 0 \\- \Sigma_A I\\0\end{bmatrix}\Sigma_A \alpha, \quad b_2 = \begin{bmatrix}0  \\ u\\q\\v\end{bmatrix}.
	\end{equation*}

	\section{Distributed Implementation}
	\label{sec:imple}
	Minimization of the convex problem in~\eqref{Final_Formulation_Z} can be carried out in multiple ways. However, we strive for a method suitable for distributed settings. Hence, we choose the \textit{Fast Iterative Shrinkage-Thresholding Algorithm (FISTA)} \cite{article:fistas}, which is an extension of the gradient method and will produce a naturally distributed solution, as detailed below. FISTA handles problems defined as 
	\begin{equation}
		\minimize_x f(x) = g(x) + h(x)
		\label{prob:FISTA}
	\end{equation}
	where $g(x) : \mathbb{R}^n \rightarrow \mathbb{R}$ is a smooth convex function, continuously differentiable with Lipschitz continuous gradient $L$, and $h(x)$ is a closed and convex function with an inexpensive proximal operator with parameter $t_\kappa$, expressed as $\textrm{prox}_{t_\kappa h}$. Under these assumptions, FISTA's solution for Problem~\eqref{prob:FISTA} is given iteratively as 
	\begin{equation}
		\begin{aligned}
			y = x^{\left(\kappa-1\right)}+\frac{\kappa-2}{\kappa+1}\Big(x^{\left(\kappa-1\right)}-x^{\left(\kappa-2\right)}\Big)\\
			x^{\left(\kappa\right)}= \textrm{prox}_{t_\kappa h}(y-t_\kappa \nabla g(y)).
		\end{aligned}
		\label{eq:fista_method}
	\end{equation} This method has a convergence rate $f(x^{\left(\kappa\right)})-f^{*}$  (where $f^{*}$ is  the optimal value) of $\mathcal{O}(1/\kappa^2)$, for a fixed step size of $t_\kappa=1/L$.

	\subsection{Instantiation}\label{sec:FISTA_instantiation}
	For Problem~\eqref{Final_Formulation_Z} to be solved with FISTA, it is necessary to define $g(x)$ and $h(x)$, along with $\textrm{prox}_{t_k h}$ and $t_k$. In this case, $g(x)$ is the quadratic cost function $\frac{1}{2}z^TMz-b^Tz$. It is noted that a function is said to have a Lipschitz Continuous Gradient with Lipschitz constant $L$ if the following condition is true
	\begin{equation}
		\label{def:Lipschitz_Continuous}
		\| \nabla f(x) - \nabla f(y) \| \leq L \|x-y\|
	\end{equation}
	and quadratic functions have a Lipschitz continuous gradient. Besides, the required gradient $\nabla g(z)$ is straightforward, given as $\nabla g(z) = Mz-b$.
	
	Term $h(x)$ corresponds to the {indicator function}\footnote{We define the indicator function as $	I_{\mathcal{Z}}(z) = \begin{cases} 0, \qquad \text{if } z \in \mathcal{Z}\\ +\infty , \quad \text{otherwise} \end{cases}$} $I_{\mathcal{Z}}(z)$ for set~$\mathcal{Z}$. {From \hbox{\cite{article:lecture_FISTA_aux}}, when} $h(x)$ is the indicator function of a closed convex set $\mathcal{Z}$, $\textrm{prox}_{t_k h}$ is the projection of $z$ on set $\mathcal{Z}$, $P_{\mathcal{Z} }(x)$, defined as 
	\begin{equation*}
		\label{def:projection}
		P_{\mathcal{Z}}(z) = \argmin_{x\in \mathcal{Z}} \| x-z\|^2,
	\end{equation*}
	which for our constraints with general structure $\{z : \| z \|\leq d \}$ is implemented as
	\begin{equation*}
		P_{\mathcal{Z}}(z) = 
		\begin{cases}z,\quad \textrm{if} \quad\|z\|\leq d \\\frac{z}{\|z\|}d,\quad \textrm{if} \quad\|z\|> d.
		\end{cases}
	\end{equation*}
	Finally, it is necessary to obtain a value for  $1/L$, where $L$ is the Lipschitz constant of function $\frac{1}{2}z^TMz-b^Tz$. Resorting to \cite{article:Soares15}, we obtain an upper-bound for $L$ as (see \ref{Appendix_Lipschitz_cte} for the complete derivation)
	\begin{equation*}
		\begin{split}
			L	&\leq \tfrac{1}{\sigma_N^2}2 \delta_{\textrm{max}}+\tfrac{1}{\sigma_A^2} \max_{i \in \mathcal{V}} |\mathcal{A}_i|+\tfrac{1}{\sigma_V^2}2\delta^V_{\textrm{max}}(T)+K
		\end{split},
	\end{equation*}
	where $\delta_{\textrm{max}}$ is the maximum node degree of the network and $ \delta^V_{\textrm{max}}$ is either $2$ for any $T>2$, $1$ for $T=2$ and $0$ for no time window ($T=1$). The term $\max_{i \in \mathcal{V}} |\mathcal{A}_i|$ corresponds to the maximum number of anchors connected to a node. Finally, $K = \tfrac{1}{\sigma_N^2} +\tfrac{1}{\sigma_A^2} +\tfrac{1}{\sigma_V^2} $, where $\sigma_N$, $\sigma_A$ and $\sigma_V$ are the maximum values of matrices $\Sigma_N$, $\Sigma_A$ and $\Sigma_V$, respectively. Distributed computation of $L$ can be attained, as suggested in~\cite{article:Soares15}, by a diffusion algorithm.
	
	{Thus, we obtain a distributed method. The explicit algorithm and the full derivation can be found in \hbox{\ref{appendix_DistributedAlgo}}. Recalling that we stay under a convex formulation, it is noted that the method} may be initialized with any feasible value without jeopardizing convergence. We further emphasize that each node on the network computes its own estimate, relying only on distances and angles taken with respect to its neighbors (nodes or anchors) and on its own velocity measurements. Notice that quantities related to other nodes are always within the set $\mathcal{A}_i$ (neighbor anchors of node $i$) or $\mathcal{E}_i$ (neighbor nodes of node $i$). At the end of each iteration (in $\kappa$), nodes also need to share their current estimates, but since this is only needed amongst neighbors, the method remains distributed.

	\section{Parameter Estimation}
	\label{sec:parameters}

Until this point, the parameters $\sigma$ and $\varkappa$ of error distributions have been assumed to be known a priori. In real applications, these parameters could be estimated from previous experiments. However, this can also be included in the algorithm so that no previous knowledge is necessary. With each new measurement, new estimates for $\varkappa$ and $\sigma$ are computed, which desirably converge to the true values. 

Maximum likelihood estimators for the parameters of each distribution are easily found in the literature. Given a random variable $X$, normally distributed with unknown mean and variance $\sigma^2$, the MLE for the latter parameter is \cite{book:statistics}
\begin{equation}
	\label{variance_estimate}
	\hat{\sigma}^2 = \frac{1}{n-1}\sum_{i=1}^{n}(X_i-\bar{X})^2.
\end{equation}
For estimation of the concentration parameter~$\varkappa$ we have \cite{article:kappa_estimate}
\begin{equation}
	\label{kappa_estimate}
	\hat{\varkappa} = \frac{\bar{\gamma}(p-\bar{\gamma}^2)}{1-\bar{\gamma}^2},
\end{equation}
where $ \bar{\gamma} = \frac{\|\sum_{i=1}^{n}x_i\|}{n} $, $x_i$ are independent and identically distributed sample unit vectors drawn from a Von Mises-Fisher distribution $vMF(\mu, \kappa)$, and $p$ is the dimension of $x_i$.

\subsection{Instantiation}

Evidently, there is no access to the true values of our measurements to accurately estimate $\varkappa$ and $\sigma$. However, after each computation of position estimates by the algorithm, it is possible to compute estimated distances, angles and velocities which can be used as "true values". So, exact values of $\varkappa$ and $\sigma$ are not to be expected, but if the position estimates are close to the true ones, then reasonable approximations should be obtained. However, {note that we rely on our position estimates to estimate the distribution parameters and subsequently rely on the estimated parameters to compute the estimated positions. So, a poor estimated value on either counterpart could affect both estimation procedures (see \hbox{Section~\ref{sec:simulation}}).}

Denoting by $\hat{x}$ the value obtained from our position estimation, it is possible to compute the resulting distances and bearings as 
\begin{equation*}
	\hat{d}_{ij}(t) = \|\hat{x}_i(t)-\hat{x}_j(t)\|, \quad \hat{u}_{ij}(t) = \frac{\hat{x}_i(t)-\hat{x}_j(t)}{\|\hat{x}_i(t)-\hat{x}_j(t)\|}.
\end{equation*}
Given the obtained measurements $d_{ij}$ and $u_{ij}$, \eqref{variance_estimate} and \eqref{kappa_estimate} {can now be used for variance and concentration parameter estimation, respectively, as}
\begin{equation*}
	\hat{\sigma}_{ij}^2(t) = {\frac{1}{t-1}}\sum_{\tau=1}^{t}(d_{ij}(\tau)-\hat{d}_{ij}(\tau))^2, \quad 	\hat{\varkappa}_{ij}(t) = \frac{\bar{\gamma}_{ij}(t)(p-\bar{\gamma}_{ij}^2(t))}{1-\bar{\gamma}_{ij}(t)^2},
\end{equation*} where 
\begin{equation}
	\bar{\gamma}_{ij}(t) = \frac{1}{t}\Bigg\lVert\sum_{\tau=1}^{t} \frac{u_{ij}(\tau)- \hat{u}_{ij}(\tau)}{\| u_{ij}(\tau)- \hat{u}_{ij}(\tau)\|}\Bigg\rVert.
\end{equation}The same reasoning can be applied to node-anchor edges, with similar resulting equations.

{Unlike \mbox{$\hat{d}_{ij}$} and \mbox{$\hat{u}_{ij}$}, estimation of velocity measurements \mbox{$\hat{\beta}_i$} requires a more careful approach. While the former requires only current estimates of the position, the latter requires positions at different time steps. Since velocity measurements can change with a high rate over time, a naive approach like averaging over a sliding window does not produce good estimates. Thus, an approach in that vein would only yield accurate estimates for linear constant trajectories.} This problem is considered in \cite{article:velocity_estimation} and the proposed solution is based on central finite differences with a time lag and rejection of high-frequency noise. The proposed estimation is then given as \begin{equation}
	\begin{aligned}
		\hat{\beta}_i(t)\Delta T = & \frac{ 5(x_i(t-3) - x_i(t-5))}{32} +\frac{ 4(x_i(t-2) - x_i(t-6))}{32}+ \frac{ (x_i(t-1) - x_i(t-7))}{32}.
	\end{aligned}
\end{equation} Recalling that $\hat{\beta}_i(t)$ is decomposed as $\hat{\beta}_i(t)=\hat{v}_i(t)\hat{V}_i(t)$, estimates for variance and concentration are given as 
\begin{equation*}
	\hat{\sigma}_{i}^2(t) ={\frac{1}{t-1}}\sum_{\tau=1}^{t}(V_i(\tau)-\hat{V}_i(\tau))^2, \quad \hat{\varkappa}_{i}(t) = \frac{\bar{\gamma}_{i}(t)(p-\bar{\gamma}_{i}^2(t))}{1-\bar{\gamma}_{i}^2(t)},
\end{equation*}
with $\bar{\gamma}_{i}(t) = \frac{\|\sum_{\tau=1}^{t}v_{i}(\tau)- \hat{v}_{i}(\tau)\|}{t}$. However, it should be noted that the introduction of a time lag has practical consequences during implementation that should be taken into account. {In particular, this means that velocity estimation is only possible after the $7$-th time step and, consequently, \mbox{$\hat{\sigma}_{i}^2$} and \mbox{$\hat{\varkappa}_{i}$} are only available after that moment. Furthermore, some time steps are necessary before achieving sufficiently accurate estimates for any of the parameters. Thus, it is still necessary to define initial values ($\hat{\sigma}^0$ and $\hat{\varkappa}^0$) for the first iterations. However, we will show in \mbox{Section~\ref{sec:simulation}} that wrong initial estimates do not increase the positioning error when performing parameter estimation.}

The addition of parameter estimation to Algorithm~\ref{alg:ALG1} produces a final distributed and parameter-free algorithm, that can be found in \ref{appendix_B}.

	\section{Simulation Results}
	\label{sec:simulation}
	{We test the performance of our proposed method in different settings for the application of AUV missions, comparing against benchmark approaches.} 
	
	\subsection{Simulations setting}\label{subsec:setting}
	
\begin{figure}[!htb]
		\centering
		\begin{subfigure}[b]{0.45\textwidth}
			\centering
			\includegraphics[width=0.8\textwidth]{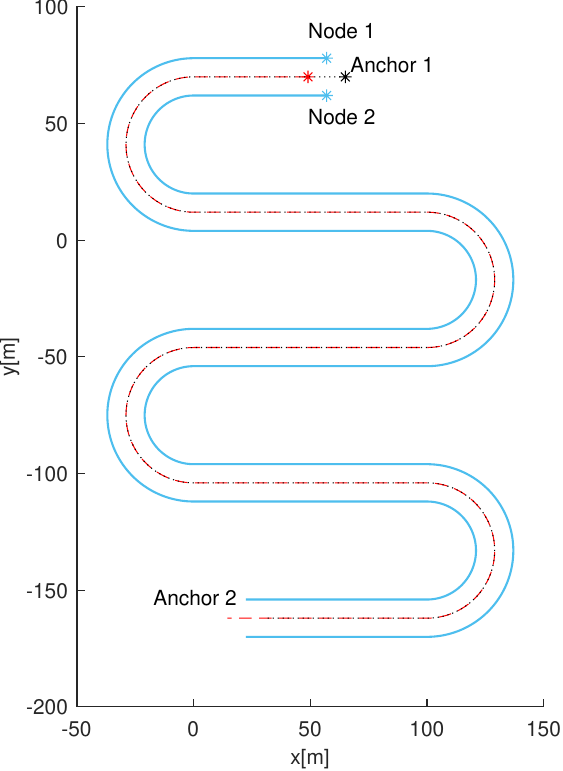}
			\caption{Lawnmower}
			\label{fig:Lawnmower_traj}
		\end{subfigure}
		\hfill
		\begin{subfigure}[b]{0.45\textwidth}
			\centering
			\includegraphics[width=\textwidth]{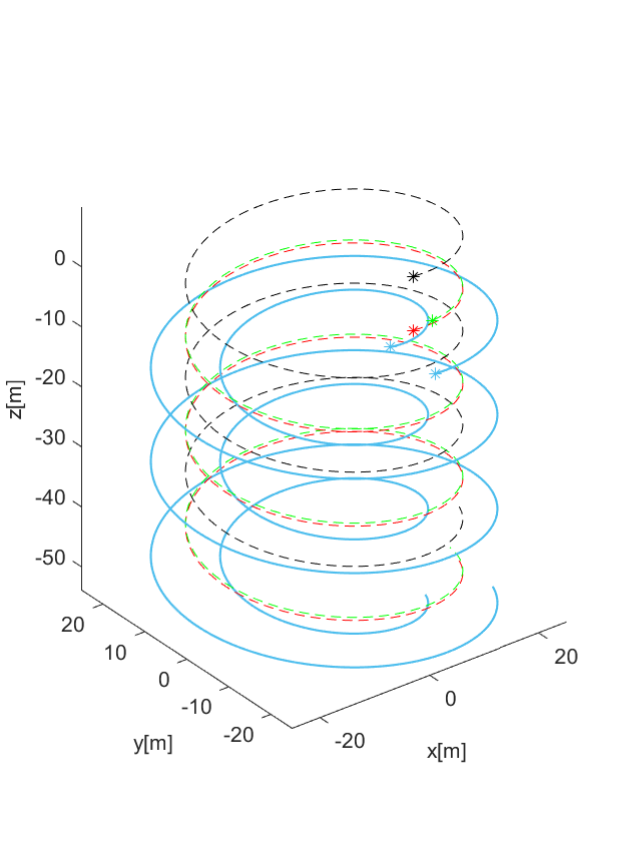}
			\caption{Helix}
			\label{fig:Helix_traj}
		\end{subfigure}
		\caption{{Lawnmower and helix trajectories. Nodes are depicted in solid blue lines and anchors in dashed lines of different colours. The start point of each trajectory is marked with a star. The lawnmower depicts two nodes, travelling on both sides of two anchors describing the same trajectory with a small lag. The helix includes two nodes, travelling on both sides of three anchors. All paths maintain the same velocity in the $x$ and $y$ axis, progressively decreasing the $z$ component.}}
		\label{fig:three graphs}
	\end{figure}

	\textbf{Trajectories.} {We evaluate the performance of the proposed algorithm over three common trajectories in AUV missions: a 2D lawnmower (\mbox{Figure~\ref{fig:Lawnmower_traj}}), a 3D helix \mbox{(Figure~\ref{fig:Helix_traj})} and a 2D lap trajectory (\mbox{Figure~\ref{fig:outlier_dist_node}}). All nodes and anchors describe the depicted trajectories at a speed of \mbox{$1m/s$} (commonly used in AUV missions \mbox{\cite{review:Wynn14,article:Fischell19}}). The three trajectories encompass different range and bearing measurements amongst themselves. Anchors are assumed to have access to their true positions, either resorting to more accurate equipment or by navigating at the surface and benefiting from GNSS.}
	
	\textbf{Noisy measurements.} {Considering noise levels obtainable with current instrumentation \mbox{\cite{review:Garcia20, article:DIESEL}}, we define \mbox{$\sigma_d = 0.5~\text{m}$} (for all distance measurements), \mbox{$\sigma_v =0.1~\text{m}/s$} (for all speed measurements), \mbox{$\kappa_a=1000$} (for all bearing measurements) and \mbox{$\kappa_v=1000$} (for all heading measurements)\footnote{Note that  \mbox{$\kappa=1000$} corresponds to approximately \mbox{$1.8\degree$}.}. Unless otherwise noted, we always consider 100 Monte Carlo trials for each simulation.}
	
	\textbf{Benchmarks.} {We want to determine two main points: whether our proposed inclusion of velocity measurements and a time window is relevant and improves the performance with respect to similar static methods; whether our resulting method is effectively competitive with current state-of-the-art solutions. To attest the former we compare with \mbox{\cite{article:Soares20}}, an hybrid convex method for ranges and bearings that does not account for velocity measurements. For the latter, we select the EKF with an equivalent data model, as it remains the state-of-the-art standard for practical AUV applications due to the trade-off between accuracy and computational power \mbox{\cite{article:Ben21,article:Zhang}}.}
	
	\textbf{Initialization and parameters.} {Parameters for the EKF are tuned by grid search for each of the considered trajectories and noise levels. EKF also requires a proper initialization to assure convergence. Considering this, the filter is initialized with the true position, contaminated with Gaussian noise of zero mean and standard deviation of $2$m \mbox{\cite{article:DIESEL}}. Our method, being convex, is independent of any initialization and this choice only affects the computation time, not the achieved accuracy (see \mbox{Section~\ref{subsec:perf}}). The only parameters in our method are the noise distribution parameters, e.g., \mbox{$\sigma$} or \mbox{$\varkappa$}. For comparison with the fine-tuned EKF, we consider the optimal case where we have access to the true parameters. When testing our parameter estimation process we consider a different setup, explained in \mbox{Section~\ref{subsec:paramest}}}.
	
	\textbf{Metrics.} In order to evaluate the localization accuracy, we use the Mean Navigation Error (MNE) defined as
	\begin{equation*}
		\text{MNE}(\hat{x},t) = \frac{1}{MN} \sum_{m=1}^{M}\sum_{i=1}^{N}  \| \hat{x}^m_i(t) -x^m_i(t) \|_2,
		\label{eq:MNE formula}
	\end{equation*}
	where $\hat{x}$ is an estimate of $x$, $M$ the number of Monte Carlo trials, $N$ the number of unknown nodes and $T$ the time interval. {Furthermore, when applicable we consider the Mean Positioning Error (MPE) as the MNE averaged over the full trajectory.}
	
	\subsection{Comparison with benchmarks} 	
	\label{subsec:bench}
	{We compare the MNE for each algorithm along the three different trajectories. The advantage introduced by the dynamic formulation with respect to the static one is clear from the decrease in MNE of \mbox{$0.2$m}. Comparison with EKF differs depending on the trajectory of the nodes. While EKF outperforms during linear segments, its error overpasses the proposed method during nonlinear motion. Mainly, our method has the advantage of keeping a constant error regardless of the vehicles' movement as it does not rely on a dynamical model of the system (unlike EKF). The corresponding MPE can be found in Table~\mbox{\ref{tab:errors}} in appendix.}
	\begin{figure}[!htb]
		\centering
		\includegraphics[width =\linewidth]{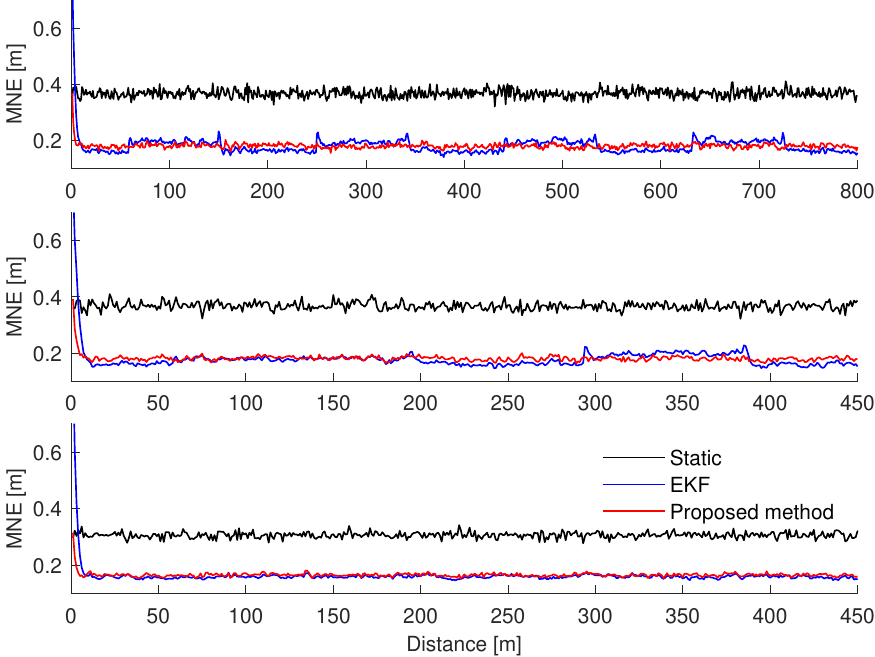}
		\caption[]{{Comparison of Mean Navigation Error between EKF, a static formulation and our algorithm, during a lawnmower trajectory (top), lap trajectory (middle) and helix trajectory (bottom), averaged over 100 MC trials. The curved segments of the trajectory are evidenced by the higher error in EKF, where our method outperforms it. For the linear parts, EKF shows better accuracy.}}
		\label{fig:Lawnmower_comparison}
	\end{figure}
	
	\subsection{Robustness to outliers}
\begin{figure}[!htb]
		\centering
		\includegraphics[width=0.9\textwidth]{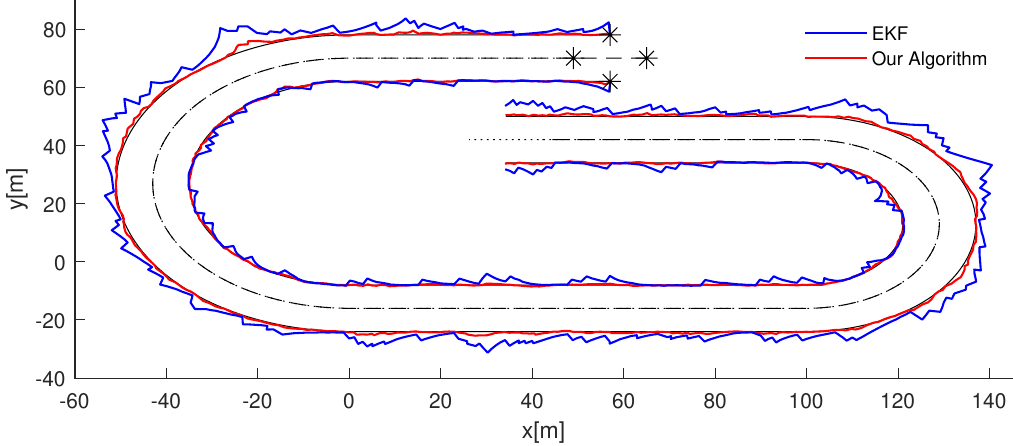}
		\caption{Outliers during a lap trajectory performed by two anchors, in dashed lines, and two nodes, in black lines, with initial positions marked as a star.  Distances between the outer node and all the other vehicles are contaminated with outliers. EKF, in blue, is affected to a higher extent and takes more time to recover than our algorithm, in red.}
		\label{fig:outlier_dist_node}
	\end{figure}

	{We consider a more realistic scenario, where outlier distance measurements with value $5d$ are introduced with a probability of $10\%$, where $d$ is the true distance.} Ranges between the outer node and all other vehicles (nodes and anchors) are subject to outliers, evidenced by the peaks in position estimates (Figure~\ref{fig:outlier_dist_node}). {The proposed algorithm presents variations of smaller magnitude when compared to EKF} and does not require as much time to recover. The difference between them also seems to accentuate during non-linear trajectory {parts with our algorithm} attaining superior robustness.
	
	\subsection{Performance of our method}\label{subsec:perf}
\begin{figure}[!htb]
		\centering
		\includegraphics[width=0.7\textwidth]{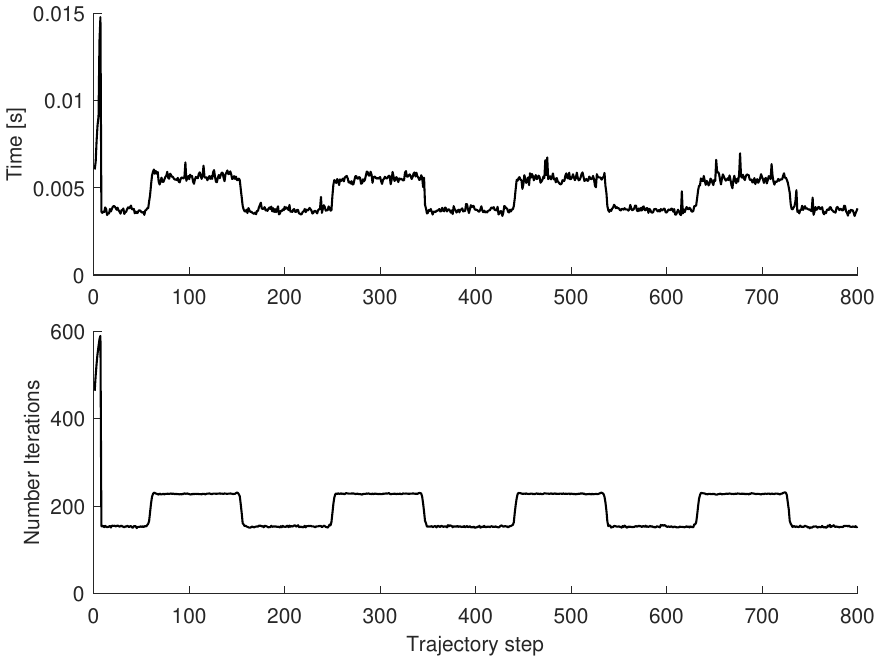}
		\caption{{Time and number of iterations spent to estimate positions, per trajectory step during a lawnmower trajectory. The values are larger during the first iterations when positions are randomly initialized and then suffer small increments during nonlinear motion.}}
		\label{fig:FistaPerf}
	\end{figure}

	{For further insight on the performance of the proposed method, we evaluate the computation time and number of iterations per step during a lawnmower trajectory \mbox{(Figure~\ref{fig:FistaPerf})}. During the first steps, while the time window is still not constant, the algorithm is initialized with a random value taking a longer time to converge. After this point, the last estimates are used for initialization and the algorithm is faster to converge (due to convexity the final estimate is independent of initialization choice and this affects only the number of iterations required for convergence). There is also a small increase in convergence time during nonlinear segments of the trajectory.}
	
	\subsection{Scalability to large networks}
\begin{figure}[!htb]
		\centering
		\includegraphics[width=0.7\textwidth]{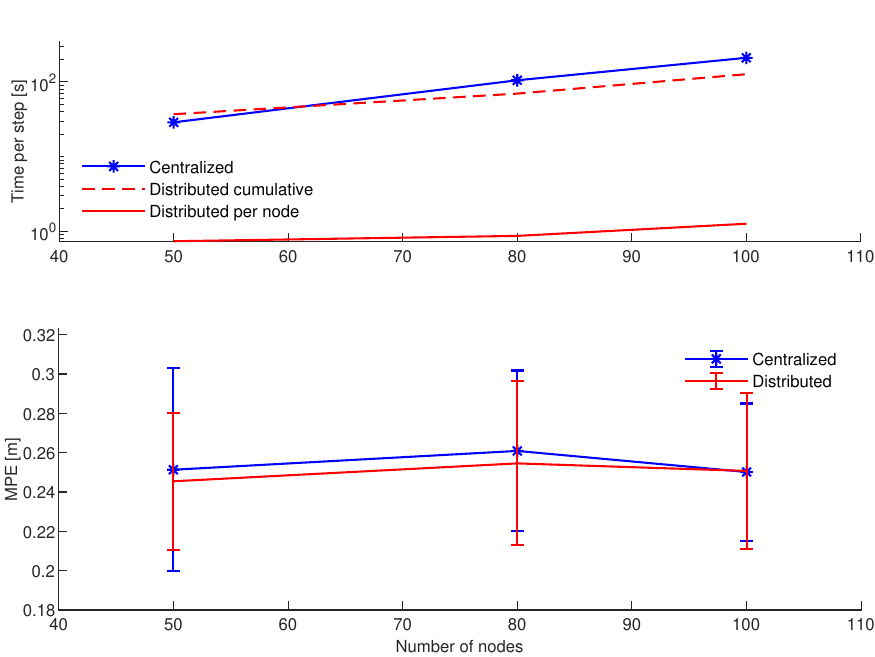}
		\caption{{Variation of the average computation time (per trajectory step) and positioning error with the number of nodes. The average node degree is kept constant over the different simulations. The blue lines refer to the centralized implementation of the proposed method and the red lines to the distributed implementation. In the top plot, the dashed line refers to the total computation time, i.e., the time for all nodes to perform estimates. The solid line is the computation time per node.}}
		\label{fig:LargeScale}
	\end{figure}

	{We test the scalability of our method by considering networks of 50, 80 and 100 nodes and comparing the centralized and distributed implementations (\mbox{Figure~\ref{fig:LargeScale}})\footnote{Details on the large scale setup can be found in \mbox{\ref{appendix:Exp}}}. The distributed implementation preserves approximately the same computation time per node, regardless of the total number of nodes of the network (if the average node degree is kept constant). Naturally, additional communication time between nodes would have to be considered in order to have an accurate estimate.}
	\subsection{Parameter estimation}\label{subsec:paramest}
\begin{figure}[!htb]
		\centering
		\includegraphics[width =0.7\linewidth]{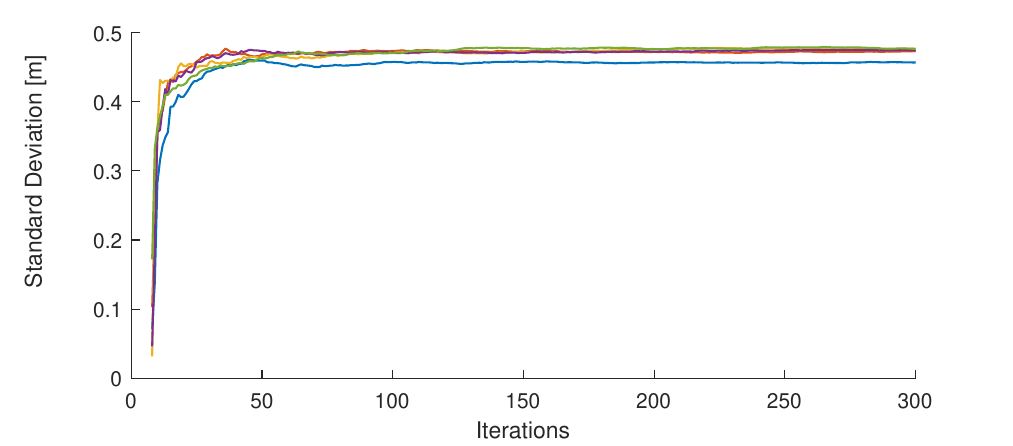}
		\caption[Estimation of $\sigma$ for distance measurements]{Evolution of estimates of $\hat{\sigma}_d$ for distance measurements. After some time, the estimates converge to values close to the true $\sigma_d = 0.5 m$. The different lines correspond to $\sigma_{ij}$ for different edges, all with the same value of STD. Note that parameter estimation does not start at iteration 0 due to the time lag in velocity estimation. }
		\label{fig:sigma}
	\end{figure}

	We start by examining the evolution of $\hat{\sigma}_d$ for all the distance measurements, when the true $\sigma_d$ was set to $0.5m$, when performing parameter estimation (Figure~\ref{fig:sigma}). The estimates converge to $\sigma_d$, as desired, but it should be noted that a considerable number of measurements is necessary before this happens. {Therefore, we require an initial value to be used in the first iterations that will be then corrected. For the next simulations we start parameter estimation at the $20$th trajectory step, using $\sigma^0$ up to that point.}
	
	\subsubsection{Impact of initial values}\label{sec:init_param}
\begin{figure}[!htb]
		\centering
		\includegraphics[width =0.8\linewidth]{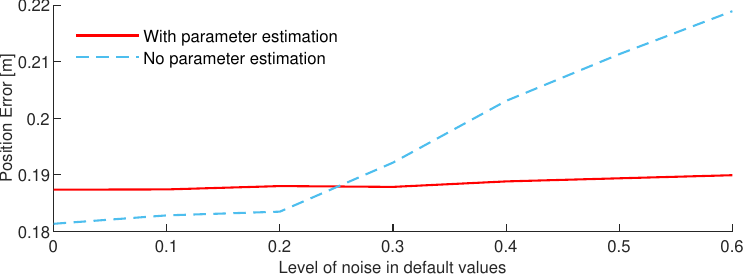}
		\caption{{Variation of MPE for different levels of noise in the initialization of all parameters. The parameters (e.g. $\sigma^0$) are initialized with the true value ($\sigma$) with additive Gaussian noise of zero mean and standard deviation as a fraction of the true value, that is $\theta \sigma$, where $\theta$ corresponds to the levels observed in the x-axis. I.e., for level $0$ we are using the true parameters as initial guesses. The algorithm with no parameter estimation uses the initial values ($\sigma^0$) throughout the entire trajectory.} }
		\label{fig:estvsnoest}
	\end{figure}

	{ Given the latter observation it is fundamental to assess the impact of the choice of the initial values, or the method does not become effectively parameter-free. We consider as initial value (e.g. $\sigma^0$) the true parameter ($\sigma$) contaminated with Gaussian noise with zero mean and increasing standard deviation, and compare the version with estimation and without (\mbox{Figure~\ref{fig:estvsnoest}}). With parameter estimation our method uses $\sigma^0$ during the first 20 iterations and then uses the estimates according to \mbox{Section~\ref{sec:parameters}}. Without estimation the method uses $\sigma^0$ throughout the entire trajectory. As expected, when $\sigma^0 \approx \sigma$, the estimation increases the error, as the estimated parameters $\hat{\sigma}$ are less accurate than  $\sigma^0$. However, as $\sigma^0$ is set further away from the true vale, the error without estimation largely increases, while the parameter estimation allows to keep the performance approximately constant. Thus, unless the distribution parameters are known with a large confidence, the estimation process is beneficial.}

	\subsubsection{Impact during trajectories} 
\begin{figure}[!htb]
		\centering
		\includegraphics[width =\linewidth]{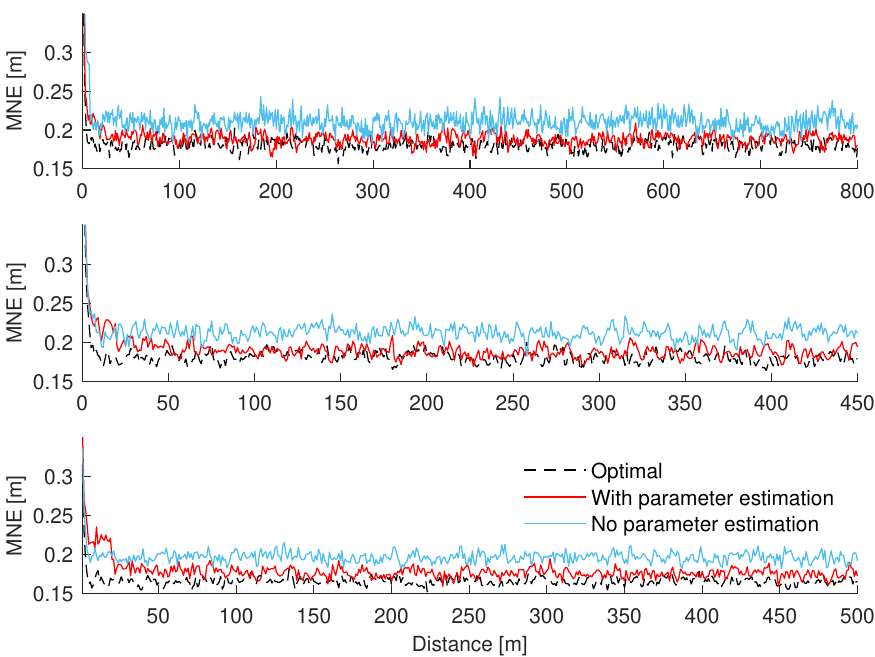}
		\caption[]{{Comparison of performance between the method with true noise parameters (\textit{Optimal}), with unknown noise parameters throughout the entire trajectory (\textit{No parameter estimation}), with initial unknown parameters and estimation procedure (\textit{With parameter estimation}). Unknown noise parameters correspond to the true ones with additive Gaussian noise with standard deviation of $0.5$ of the true value. The results correspond to a lawnmower trajectory (top), lap trajectory (middle) and helix trajectory (bottom).} }
		\label{fig:ssigma}
	\end{figure}

{Finally, we compare the algorithm with and without parameter estimation at level $0.5$ of noise in initial parameters (as defined in \mbox{Section~\ref{sec:init_param}}) and against the optimal setting where the true noise parameters are known (\mbox{Figure~\ref{fig:ssigma}}). During the initial iterations, the noisy parameters are being used in both versions, leading to a higher MNE. After this point, while algorithm with no estimation keeps a similar error, the estimation process allows to achieve an error close to the optimal one. The corresponding MPE can be found in Table~\mbox{\ref{tab:errors}} in \mbox{\ref{appendix:Exp}}.}

	\section{Concluding Remarks}
	\label{sec:concl}	
We extend the static MLE formulation for a generic mobile network localization problem from hybrid (range and bearing) measurements to the dynamic setting, in an horizon-based approach. This entails only the additional assumption that each node can measure its own velocity. We show that this extended formulation is still amenable to a tight convex relaxation and can be implemented in a distributed manner.

Comparison with EKF, a centralized state-of-the-art method, shows the superior performance of our method for cases with variable motion patterns, where for mostly linear trajectories the EKF still slightly outperforms our distributed estimator. Nonetheless, unlike EKF, our distributed convex method is proven to converge regardless of initialization and offer nearly constant accuracy throughout the trajectory, thus being more predictable than EKF. More importantly, our algorithm shows more robustness in the presence of outliers in measurements, which is a known issue in real-life scenarios.

Besides, with the proposed parameter estimation process, no parameter tuning is required, so the algorithm can be easily applied to new environments or technologies. However, it should be noted that this is a sensitive step where abnormal values could start a chain of wrong position and parameter estimates. While a formal proof would be desirable in future approaches, our simulations have consistently demonstrated that the method maintains its performance, regardless of the initialization chosen and despite the estimated parameters. Therefore, we have strong empirical evidence of robustness for the proposed method, although additional research is needed to prove its validity formally.

	\section{Acknowledgments}
	This work was supported by Recovery and Resilience Plan and NextGeneration EU Funds through Project \textit{Artificial Intelligence Fights Space Debris} [C626449889-0046305]; by LARSyS - FCT Project [UIDB/50009/2020]; by NOVA LINCS [Grant UIDB/04516/2020]; and by national funds through the FCT – Fundação para a Ciência e a Tecnologia, I.P., under the scope of the projects UIDB/00297/2020 (https://doi.org/10.54499/UIDB/00297/2020) and UIDP/00297/2020 (https://doi.org/10.54499/UIDP/00297/2020) (Center for Mathematics and Applications).

	\bibliographystyle{unsrt}
	\bibliography{Ref_db}  

\begin{thebibliography}{10}

\bibitem{article:Miller10}
P.~A. Miller, J.~A. Farrell, Y.~Zhao, and V.~Djapic.
\newblock Autonomous underwater vehicle navigation.
\newblock {\em IEEE Journal of Ocean. Eng.}, 35(3):663--678, July 2010.

\bibitem{article:Yoo14}
T.~Yoo.
\newblock {DVL}/{RPM} based velocity filter aiding in the underwater vehicle
  integrated inertial navigation system.
\newblock {\em Journal of Sensor Technology}, 4, 2014.

\bibitem{article:UAV_coop}
Lang Ruan, Guangxia Li, Weiheng Dai, Shiwei Tian, Guangteng Fan, Jian Wang, and
  Xiaoqi Dai.
\newblock Cooperative relative localization for {UAV} swarm in {GNSS}-denied
  environment: A coalition formation game approach.
\newblock {\em IEEE Internet of Things J.}, 9(13):11560--11577, 2022.

\bibitem{article:coop_indoors}
Lukas Wielandner, Erik Leitinger, and Klaus Witrisal.
\newblock {RSS-based} cooperative localization and orientation estimation
  exploiting antenna directivity.
\newblock {\em IEEE Access}, PP:1--1, 04 2021.

\bibitem{article:Pinheiro17}
Breno~C. Pinheiro, Ubirajara~F. Moreno, João T.~B. de~Sousa, and Orlando~C.
  Rodríguez.
\newblock Kernel-function-based models for acoustic localization of underwater
  vehicles.
\newblock {\em IEEE Journal of Ocean. Eng.}, 42(3):603--618, July 2017.

\bibitem{article:Dong12}
Liang Dong.
\newblock Cooperative localization and tracking of mobile ad hoc networks.
\newblock {\em IEEE Trans. Signal Process.}, 60(7):3907--3913, 2012.

\bibitem{article:EKF_2023}
Guangrun Sheng, Xixiang Liu, Yehua Sheng, Xiangzhi Cheng, and Hao Luo.
\newblock Cooperative navigation algorithm of extended {Kalman} filter based on
  combined observation for {AUVs}.
\newblock {\em Remote Sensing}, 15(2), 2023.

\bibitem{book:Kay97}
S.~M. Kay.
\newblock {\em Fundamentals of Statistical Signal Processing: Estimation
  Theory}.
\newblock Prentice Hall, 1997.

\bibitem{article:calafiore10}
G.~C. Calafiore, L.~Carlone, and M.~Wei.
\newblock Distributed optimization techniques for range localization in
  networked systems.
\newblock In {\em 49th IEEE Conf. on Decision and Control (CDC)}, pages
  2221--2226, Dec 2010.

\bibitem{article:distributed_nonconvex_costa2006}
Jose~A. Costa, Neal Patwari, and Alfred~O. Hero.
\newblock Distributed weighted-multidimensional scaling for node localization
  in sensor networks.
\newblock {\em ACM Trans. on Sens. Netw.}, 2(1):39–64, Feb 2006.

\bibitem{article:Erseghe2015}
Tomaso Erseghe.
\newblock A distributed and maximum-likelihood sensor network localization
  algorithm based upon a nonconvex problem formulation.
\newblock {\em IEEE Trans. on Signal and Information Process. over Netw.},
  1(4):247--258, 2015.

\bibitem{article:Biswas06}
Pratik Biswas and Yinyu Ye.
\newblock A distributed method for solving semidefinite programs arising from
  ad hoc wireless sensor network localization.
\newblock {\em Multiscale optimization methods and applications}, pages 69--84,
  2006.

\bibitem{YANG2022108504}
Ge~Yang et~al.
\newblock Improved robust {TOA}-based source localization with individual
  constraint of sensor location uncertainty.
\newblock {\em Signal Processing}, 196:108504, 2022.

\bibitem{article:SOCP_Tseng_2006}
P.~Tseng.
\newblock Second order cone programming relaxation of sensor network
  localization.
\newblock {\em SIAM J. on Optim.}, 18(1):156--185, 2007.

\bibitem{article:SOS_2006}
Jiawang Nie.
\newblock Sum of squares method for sensor network localization.
\newblock {\em Computational Optimization and Applications}, 43, Jun 2006.

\bibitem{article:distributed_SOCP}
Seshan Srirangarajan, Ahmed Tewfik, and Zhi-Quan Luo.
\newblock Distributed sensor network localization using {SOCP} relaxation.
\newblock {\em IEEE Trans. on Wirel. Commun.}, 7:4886--4895, 2009.

\bibitem{article:Shi10}
Q.~Shi, C.~{He}, H.~Chen, and L.~{Jiang}.
\newblock Distributed wireless sensor network localization via sequential
  greedy optimization algorithm.
\newblock {\em IEEE Trans. Signal Process.}, 58(6):3328--3340, June 2010.

\bibitem{SOARES2021}
Cláudia Soares and João Gomes.
\newblock {STRONG}: Synchronous and asynchronous robust network localization,
  under non-gaussian noise.
\newblock {\em Signal Process.}, 185:108066, 2021.

\bibitem{article:Wang23}
Zhuping Wang, Yanhao Chang, Hao Zhang, and Huaicheng Yan.
\newblock Bearing-only distributed localization for multi-agent systems with
  complex coordinates.
\newblock {\em Information Sciences}, 626:837--850, 2023.

\bibitem{article:Crouse13}
David Crouse, Richard Osborne, Krishna Pattipati, P.~Willett, and Yaakov
  Bar-Shalom.
\newblock Efficient {2D} sensor location estimation using targets of
  opportunity.
\newblock {\em Jounal of Advances in Information Fusion}, 8, Jun 2013.

\bibitem{article:convex_optimization_uav}
Liu Ming-Yong, Li~Wen-Bai, and Pei Xuan.
\newblock Convex optimization algorithms for cooperative localization in
  autonomous underwater vehicles.
\newblock {\em Acta Automatica Sinica}, 36(5):704--710, 2010.

\bibitem{article:FLORIS}
Beatriz~Quintino Ferreira, João Gomes, Cláudia Soares, and João~P. Costeira.
\newblock {FLORIS and CLORIS}: Hybrid source and network localization based on
  ranges and video.
\newblock {\em Signal Process.}, 153:355--367, 2018.

\bibitem{article:Biswas5}
P.~Biswas, H.~Aghajan, and Y.~Ye.
\newblock Semidefinite programming algorithms for sensor network localization
  using angle information.
\newblock In {\em 39th Asilomar Conference on Signals, Systems and Computers},
  pages 220--224, 2005.

\bibitem{article:SDP_benchmark}
H.~Naseri and V.~Koivunen.
\newblock Convex relaxation for maximum-likelihood network localization using
  distance and direction data.
\newblock In {\em 2018 IEEE 19th International Workshop on Signal Processing
  Advances in Wireless Communications (SPAWC)}, pages 1--5, June 2018.

\bibitem{article:Hybrid_Eren2011}
Tolga Eren.
\newblock Cooperative localization in wireless ad hoc and sensor networks using
  hybrid distance and bearing (angle of arrival) measurements.
\newblock {\em EURASIP Journal on Wirel. Commun. and Networking}, 2011, 2011.

\bibitem{article:mixed_Lin}
Zhiyun Lin, Tingrui Han, Ronghao Zheng, and Changbin Yu.
\newblock Distributed localization with mixed measurements under switching
  topologies.
\newblock {\em Automatica}, 76:251--257, 2017.

\bibitem{article:HybridRSS_Ding2021}
Weizhong Ding, Shengming Chang, and Jun Li.
\newblock A novel weighted localization method in wireless sensor networks
  based on hybrid {RSS/AoA }measurements.
\newblock {\em IEEE Access}, PP:1--1, 2021.

\bibitem{article:hybrid_RSSAoA_Tomic2019}
Slavisa Tomic, Marko Beko, and Milan Tuba.
\newblock A linear estimator for network localization using integrated {RSS}
  and {AOA} measurements.
\newblock {\em IEEE Signal Process. Lett.}, 26(3):405--409, 2019.

\bibitem{article:Soares20}
Cláudia Soares, Filipa Valdeira, and João Gomes.
\newblock Range and bearing data fusion for precise convex network
  localization.
\newblock {\em IEEE Signal Process. Lett.}, 27:670--674, 2020.

\bibitem{article:Kumar16}
Sandeep Kumar, Raju Kumar, and Ketan Rajawat.
\newblock Cooperative localization of mobile networks via velocity-assisted
  multidimensional scaling.
\newblock {\em IEEE Trans. Signal Process.}, 64(7):1744--1758, 2016.

\bibitem{article:RangeVelRelative}
Yingsheng Fan, Kai Ding, Xiaogang Qi, and Lifang Liu.
\newblock Cooperative localization of {3D} mobile networks via relative
  distance and velocity measurement.
\newblock {\em IEEE Communications Lett.}, 25(9):2899--2903, 2021.

\bibitem{article:distributed_DILAND}
Usman Khan, Karm Soummya, and Jose Moura.
\newblock {DILAND}: An algorithm for distributed sensor network localization
  with noisy distance measurements.
\newblock {\em IEEE Trans. Signal Process.}, 58:1940--1947, 2010.

\bibitem{article:Simonetto}
Andrea Simonetto and Geert Leus.
\newblock Distributed maximum likelihood sensor network localization.
\newblock {\em IEEE Trans. Signal Process.}, 62(6):1424--1437, 2014.

\bibitem{article:DistribCoop_Salari}
Soheil Salari, Il-Min Kim, and Francois Chan.
\newblock Distributed cooperative localization for mobile wireless sensor
  networks.
\newblock {\em IEEE Wirel. Commun. Lett.}, 7(1):18--21, 2018.

\bibitem{article:fistas}
A.~Beck and M.~Teboulle.
\newblock A fast iterative shrinkage-thresholding algorithm with application to
  wavelet-based image deblurring.
\newblock In {\em 2009 IEEE International Conference on Acoustics, Speech and
  Signal Processing}, pages 693--696, April 2009.

\bibitem{article:lecture_FISTA_aux}
Neal Parikh and Stephen Boyd.
\newblock Proximal algorithms.
\newblock {\em Foundations and trends in Optimization}, 1(3):127--239, 2014.

\bibitem{article:Soares15}
Cláudia Soares, J.~Xavier, and J.~Gomes.
\newblock Simple and fast convex relaxation method for cooperative localization
  in sensor networks using range measurements.
\newblock {\em IEEE Trans. Signal Process.}, 63(17):4532--4543, Sept 2015.

\bibitem{book:statistics}
Douglas~C. Montgomery and George~C. Runger.
\newblock {\em Applied Statistics and Probability for Engineers}.
\newblock John Wiley and Sons, Inc., 2003.

\bibitem{article:kappa_estimate}
Arindam Banerjee, Inderjit~S. Dhillon, Joydeep Ghosh, and Suvrit Sra.
\newblock Clustering on the unit hypersphere using von {Mises-Fisher}
  distributions.
\newblock {\em J. Mach. Learn. Res.}, 6:1345--1382, 2005.

\bibitem{article:velocity_estimation}
C.~Soares, J.~Gomes, B.~Q. Ferreira, and J.~P. Costeira.
\newblock {LocDyn}: Robust distributed localization for mobile underwater
  networks.
\newblock {\em IEEE Journal of Ocean. Eng.}, 42(4):1063--1074, Oct 2017.

\bibitem{review:Wynn14}
Russell~B. Wynn and et~al.
\newblock Autonomous underwater vehicles ({AUVs}): Their past, present and
  future contributions to the advancement of marine geoscience.
\newblock {\em Marine Geology}, 352:451--468, 2014.
\newblock 50th Anniversary Special Issue.

\bibitem{article:Fischell19}
Erin~Marie Fischell, Nicholas~Rahardiyan Rypkema, and Henrik Schmidt.
\newblock Relative autonomy and navigation for command and control of low-cost
  autonomous underwater vehicles.
\newblock {\em IEEE Robotics and Automation Letters}, 4(2):1800--1806, 2019.

\bibitem{review:Garcia20}
et~al. Gonz{\'a}lez-Garc{\'\i}a.
\newblock Autonomous underwater vehicles: Localization, navigation, and
  communication for collaborative missions.
\newblock {\em Applied sciences}, 10(4):1256, 2020.

\bibitem{article:DIESEL}
Cl{\'{a}}udia Soares, Pusheng Ji, Jo{\~{a}}o~Pedro Gomes, and Ant{\'{o}}nio
  Pascoal.
\newblock {DIeSEL}: Distributed self-localization of a network of underwater
  vehicles, 2017.

\bibitem{article:Ben21}
Yueyang Ben, Yan Sun, Qian Li, and Xinle Zang.
\newblock A novel cooperative navigation algorithm based on factor graph with
  cycles for auvs.
\newblock {\em Ocean Engineering}, 241:110024, 2021.

\bibitem{article:Zhang}
Lingling Zhang, Shijiao Wu, and Chengkai Tang.
\newblock Cooperative positioning of underwater unmanned vehicle clusters based
  on factor graphs.
\newblock {\em Ocean Engineering}, 287:115854, 2023.

\bibitem{article:laplacian_eigenvalues}
William N.~Anderson Jr. and Thomas~D. Morley.
\newblock Eigenvalues of the {Laplacian} of a graph.
\newblock {\em Linear and Multilinear Algebra}, 18(2):141--145, 1985.

\end{thebibliography}
	\newpage
	\appendix

\counterwithin*{equation}{section}
\renewcommand{\theequation}{\thesection.\arabic{equation}}

\section{Upper bound on Lipschitz Constant}
\label{Appendix_Lipschitz_cte}

First, the gradient difference is manipulated until an expression with the structure of \eqref{def:Lipschitz_Continuous} is obtained as
\begin{equation}
	\begin{split}
		\| \nabla f(x) - \nabla f(y) \|& = \|( Mx -b) - (My -b) \| \\
		&=\| M (x -y) \| \\
		&\leq \| M  \|\|(x -y) \| \quad .
	\end{split}
	\label{L_for_quadratic}
\end{equation}

The last step of \eqref{L_for_quadratic} corresponds to the definition of Lipschitz continuity with $L=\|M\|$, where $\|M\|$ is the matrix norm of $M$ defined as 
\begin{equation*}
	\|M\| = \sup \Bigg\{\frac{\|Mx\|}{\|x\|}: x \in R^n, |x| \neq 0\Bigg\} \quad.
\end{equation*}

Further manipulations are necessary to obtain a bound for this constant; these are shown in \eqref{Lip_manipulation} and the justification is presented below.

\begin{equation}
	\label{Lip_manipulation}
	\begin{split}
		L=&\|M\|\\
		\stackrel{\textrm{1}}{\leq}& \|M_1\|+\|M_2\|+\|M_3\|\\
		\stackrel{\textrm{2}}{=}&\sqrt{\lambda_{\textrm{max}}(M_1^TM_1)}+\sqrt{\lambda_{\textrm{max}}(M_2^TM_2)}+\sqrt{\lambda_{\textrm{max}}(M_3^TM_3)}\\ 
		\stackrel{\textrm{3}}{=}&\lambda_{\textrm{max}}(M_1)+\lambda_{\textrm{max}}(M_2)+\lambda_{\textrm{max}}(M_3)\\
		\stackrel{\textrm{4}}{=}&\lambda_{\textrm{max}}(\Sigma_N^2 AA^T +\Sigma_N^2 ) + \lambda_{\textrm{max}}(\Sigma_A^2 EE^T +\Sigma_A^2 )\\&+ \lambda_{\textrm{max}}(\Sigma_V^2 NN^T +\Sigma_V^2 ) \\
		\stackrel{\textrm{5}}{\leq}&\lambda_{\textrm{max}}(\Sigma_N^2 AA^T) + \lambda_{\textrm{max}}(\Sigma_A^2 EE^T )+ \lambda_{\textrm{max}}(\Sigma_V^2 NN^T ) \\&+  \lambda_{\textrm{max}}(\Sigma_N^2)+  \lambda_{\textrm{max}}(\Sigma_A^2)+  \lambda_{\textrm{max}}(\Sigma_V^2) \\
		\stackrel{\textrm{6}}{\leq}&\tfrac{1}{\sigma_N^2}\lambda_{\textrm{max}}( AA^T) + \tfrac{1}{\sigma_A^2}\lambda_{\textrm{max}}(EE^T)+ \tfrac{1}{\sigma_V^2}\lambda_{\textrm{max}}( NN^T)\\&+(\tfrac{1}{\sigma_N^2} +\tfrac{1}{\sigma_A^2} + \tfrac{1}{\sigma_V^2} ) \\ 
	\end{split}
\end{equation}

\begin{itemize}
	\item[]$\stackrel{\textrm{1}}{}$ For 2 matrices $A$ and $B$, $\|A+B\| \leq\|A\| +\|B\|$ (triangle inequality)
	\item[]$\stackrel{\textrm{2}}{}$ $\|A\|=\sqrt{\lambda_{\textrm{max}}(A^TA)}$
	\item[]$\stackrel{\textrm{3}}{}$ $M_1$,$M_2$ and $M_3$ are symmetric ($A^T=A$), thus verifying  $\sqrt{\lambda_{\textrm{max}}(A^TA)}=\lambda_{\textrm{max}}(A)$
	\item[]$\stackrel{\textrm{4}}{}$ $\lambda_{\textrm{max}}(A^TB) = \lambda_{\textrm{max}}(AB^T)$. For example, \begin{equation*}\begin{aligned}
			\lambda_{\textrm{max}}( M_1) &=\lambda_{\textrm{max}} \Bigg(	\begin{bmatrix}A^T \\- I  \\ 0\\0\end{bmatrix}
			\begin{bmatrix}\Sigma_N^2 A & - \Sigma_N^2  & 0& 0\end{bmatrix} \Bigg) \\& =\lambda_{\textrm{max}} \Bigg(	
			\begin{bmatrix}\Sigma_N^2 A & - \Sigma_N^2  & 0& 0\end{bmatrix}\begin{bmatrix}A^T \\- I  \\ 0\\0\end{bmatrix}\Bigg)
	\end{aligned}\end{equation*}
	\item[]$\stackrel{\textrm{5}}{}$ Weyl's inequality
	\item[]$\stackrel{\textrm{6}}{}$ Taking $\tfrac{1}{\sigma_N^2}$, $\tfrac{1}{\sigma_A^2}$,$\tfrac{1}{\sigma_V^2}$ as the maximum values of diagonal matrices $\Sigma_N$, $\Sigma_A$ and $\Sigma_V$.
\end{itemize}

The last step is intentionally separated, since it relates with particular characteristics of the network problem. First, the \textit{Laplacian matrix} $L(\mathcal{G})$ of a graph $\mathcal{G}$ with incidence matrix $C$ is defined as $L (\mathcal{G}) = CC^T$. Recalling that matrix $A$ is in fact an arc-node incidence matrix of a virtual network, it is concluded that $AA^T$ is the Laplacian matrix of such network. Additionally, it was proven in \cite{article:laplacian_eigenvalues} that the maximum eigenvalue of $L (\mathcal{G})$ is upper-bounded by twice the maximum node degree of $\mathcal{G}$. Combining these two considerations, it follows that $\lambda_{\textrm{max}}( AA^T)\leq2 \delta^V_{\textrm{max}}$, where $ \delta^V_{\textrm{max}}$ denotes the maximum node degree of the virtual network. Finally, under the assumption that the set of edges remains constant under the time window, $\delta^V_{\textrm{max}}$ is equivalent to $\delta_{\textrm{max}}$, the maximum node degree of the true network.

The same considerations are applicable to $\lambda_{\textrm{max}}( NN^T)$, with $\lambda_{\textrm{max}}( NN^T)\leq2 \delta^V_{\textrm{max}}$. Here, $ \delta^V_{\textrm{max}}$ is either $2$ for any $T>2$, $1$ for $T=2$, and $0$ for no time window ($T=1$). It will then be designated as $\delta^1_{\textrm{max}}(T)$, a function of time window. Finally, $\lambda_{\textrm{max}}(EE^T)$ is noted to be upper bounded by the maximum cardinality over all the sets $\mathcal{A}_i$ \cite{article:Soares15}. In other words, this value corresponds to the maximum number of anchors connected to a node.

With $K = \tfrac{1}{\sigma_N^2} +\tfrac{1}{\sigma_A^2} +\tfrac{1}{\sigma_V^2} $ and the above considerations, the previous upper bound may transformed as
\begin{equation*}
	\begin{split}
		L	&\leq \tfrac{1}{\sigma_N^2}\lambda_{\textrm{max}}( AA^T) + \tfrac{1}{\sigma_A^2}\lambda_{\textrm{max}}(EE^T)   +\tfrac{1}{\sigma_V^2}\lambda_{\textrm{max}}( NN^T)+K \\
		&\leq \tfrac{1}{\sigma_N^2}2 \delta_{\textrm{max}}+\tfrac{1}{\sigma_A^2} \max_{i \in \mathcal{V}} |\mathcal{A}_i|+\tfrac{1}{\sigma_V^2}2\delta^1_{\textrm{max}}(T)+K\quad.
	\end{split}
\end{equation*}

\section{Distributed algorithm}
\label{appendix_DistributedAlgo}

\subsection{Distributed Algorithm}
We present the algorithm for the algorithm for distributed method in Algorithm~\ref{alg:ALG1}. Note that the algorithm contains an outer iteration over time $t$ and an inner iteration $\kappa$ for the solution of problem~\eqref{Final_Formulation_Z} with FISTA. Finally, iterations in $\tau$ correspond to the time window. We define $C_{(i\sim j,i)}$ as the element of $C$ on row referring to edge $i\sim j$ and column referring to node $i$; $N_i$ as the sub-matrix of $N$ whose rows correspond to node $i$; $\textbf{s}_i$ and $\textbf{x}_i$ as the concatenation of $s_i(\tau)$ and $x_i(\tau)$ over the time window; $\textbf{w}_i$ and $\textbf{y}_i$ as the concatenation of $w_{ik}(\tau)$ and $y_{ij}(\tau)$ over the respective edges and time window. Besides, the projections for each component are defined as
$\mathcal{Y}_{ij}(t) = \{y_{ij}(t) : \| y(t) \|\leq d_{ij}(t) \}$,
$\mathcal{W}_{ik}(t) = \{w_{ik}(t) : \| w_{ik}(t) \|\leq r_{ik}(t) \}$
and $\mathcal{S}_{i}(t) = \{s_i(t) : \| s_i(t) \|\leq v_{i}(t)
\}$. The elements $F_1(\tau)$ through $F_4(\tau)$ for the update of $x_i^{\kappa+1}(\tau)$ are given as
\begin{equation*}
	F_1(\tau) =  \hat{x}_i^\kappa(\tau) \Big(L - \sum_{i\sim j \in \mathcal{E}_i} \frac{1}{\sigma_{ij}^2} -\sum_{k \in \mathcal{A}_i}\frac{1}{\varsigma_{ik}^2}\Big),
\end{equation*}
\begin{equation*}
	F_2(\tau) = \sum_{i\sim j \in \mathcal{E}_i} \frac{1}{\sigma_{ij}^2}(\hat{x}_j^{\kappa}(\tau)+C_{(i\sim j,i)}(\tau)\hat{y}_{ij}^\kappa(\tau)),
\end{equation*}
\begin{equation*}
	F_3(\tau) = \sum_{k \in \mathcal{A}_i}\frac{1}{\varsigma_{ik}^2}(\hat{w}_{ik}^\kappa(\tau)+\alpha_k(\tau))
\end{equation*}
and
\begin{equation*}
	F_4(\tau)=\big[N_i^T\Sigma_V^2(\hat{\textbf{s}}_i^{\kappa} - N_i\hat{\textbf{x}}_i^\kappa) \big]_{\tau},
\end{equation*} where $\big[v\big]_\tau$ is the element of a vector $v$, relative to instant $\tau$. The update depends on the previous iteration $\kappa$ and is split into four terms for convenience of the reader: $F_1(\tau)$ relates to the node position; $F_2(\tau)$ and $F_3(\tau)$ relate to distance edges between the node and its neighbor nodes or anchors (respectively); and $F_4(\tau)$ to the node velocity.

\begin{algorithm}[H]	\caption{Hybrid Convex Localization}\label{alg:ALG1}
	\begin{algorithmic}
		\REQUIRE $L$; 
		\STATE $\{d_{ij}(t), u_{ij}(t):i \sim j \in \mathcal{E}\}$; $\{r_{ik}(t), q_{ik}(t):i \in \mathcal{V},k \in \mathcal{A}\}$; $\{v_i(t) \in \mathcal{V} \} $
		\ENSURE{$\hat{x}=\{\hat{x}_i\}$ }
		\STATE Choose $z^0=(x^0,y^0,w^0,s^0)=0$;
		\FORALL{t}
		\STATE{$\kappa=1$; $z^1=z^0$}
		\WHILE{some stopping criterion is not met, each node $i$}
		\STATE$\hat{z}_i^{\kappa}=z_i^\kappa+\frac{\kappa-1}{\kappa}(z_i^\kappa-z_i^{\kappa-1})$
		\FORALL{$t-T_0\leq \tau \leq t$}
		\STATE{$x_i^{\kappa+1}(\tau) = \frac{1}{L}\Big[ F_1(\tau)+F_2(\tau)+F_3(\tau)+F_4(\tau)  \Big] $}
		\ENDFOR			
		\FORALL{$i\sim j \in \mathcal{E}_i$ and $t-T_0\leq \tau \leq t$}
		\STATE{$y_{ij}^{\kappa+1}(\tau) = P_{\mathcal{Y}_{ij}(\tau)}\Big( \frac{L-1/\sigma_{ij}^2}{L}\hat{y}_{ij}^{\kappa}(\tau) + \frac{\tilde{u}_{ij}(\tau)}{L} + \frac{1}{L\sigma_{ij}^2}C_{(i\sim j,i)}(\tau)(\hat{x}_i^\kappa(\tau)-\hat{x}_j^\kappa(\tau))\Big)$}		
		\ENDFOR
		\FORALL{$k \in \mathcal{A}_i$ and $t-T_0\leq \tau \leq t$}
		\STATE{$w_{ik}^{\kappa+1}(\tau) = P_{\mathcal{W}_{ik}(\tau)}\Big( \frac{L-1/\varsigma_{ik}^2}{L}\hat{w}_{ik}^{\kappa}(\tau)+\frac{1}{L\varsigma_{ik}^2}(\hat{x}_i^\kappa(\tau)-a_k(\tau))+\frac{\tilde{q}_{ik}(\tau)}{L}\Big)$}		
		\ENDFOR
		\IF{$T_0>1$}
		\FORALL{$t-T_0+1\leq \tau \leq t$}
		\STATE{$s_{i}^{\kappa+1}(\tau) = P_{\mathcal{S}_{i}}\Big(\frac{L-1/\sigma_{ij}^2}{L}+\frac{1}{L\sigma_{i}^2}(\hat{x}_i^\kappa(\tau)-\hat{x}_i^\kappa(\tau-1)) +\frac{\tilde{v}_{i}(\tau)}{L} \Big)$}
		\ENDFOR
		\ENDIF
		\STATE $z_i^{\kappa+1}=(\textbf{x}_i^{\kappa+1},\textbf{y}_i^{\kappa+1},\textbf{w}_i^{\kappa+1},\textbf{s}_i^{\kappa+1})$
		\STATE $\kappa=\kappa+1$
		\STATE broadcast $x_i$ to all neighbors
		\ENDWHILE
		\STATE \RETURN $\hat{x} = \{\textbf{x}_i^{\kappa+1}\}$
		\ENDFOR
	\end{algorithmic}
\end{algorithm}

\subsection{Derivation of the distributed algorithm}
In this section, we present the full derivation of Algorithm~\ref{alg:ALG1}. We start by replacing the gradient $\nabla g(z) = Mz-b$ in the update equations for FISTA \eqref{eq:fista_method}, obtaining the two updates as
\begin{equation}
	\begin{aligned}
		\hat{z}^{\kappa} =& z^{\kappa}+\frac{\kappa-1}{\kappa}\Big(z^{\kappa}-z^{\kappa-1}\Big)\\
		z^{\kappa+1} =& P_{\mathcal{Z}} \Big(\hat{z}^\kappa-\frac{1}{L}(M\hat{z}^\kappa-b)\Big).
	\end{aligned}
\end{equation} We shall further develop the second equation, in order to retrieve the update for each node $i$. Considering the update for each component of $z$ and the definitions of matrices $M$ and $b$, we get
\begin{equation}
	\begin{aligned}
		x^{\kappa+1} &=  \Big[I-\frac{1}{L}(A^T\Sigma_N^2A+E^T\Sigma_A^2E+N^T\Sigma_V^2N)\Big]\hat{x}^{\kappa} +\frac{1}{L}\Sigma_N^2A^T\hat{y}^{\kappa}+\frac{1}{L}\Sigma_A^2E^T(\hat{w}^{\kappa}+\alpha) + \frac{1}{L}\Sigma_V^2N^T\hat{s}^{\kappa}  \\
		y^{\kappa+1} &= P_{\mathcal{Y}} \Big((I-\frac{1}{L}\Sigma_N^2)\hat{y}^\kappa+\frac{1}{L}\Sigma_N^2A\hat{x}^\kappa + \frac{1}{L}u \Big)\\
		w^{\kappa+1} &= P_{\mathcal{W}} \Big( (I-\frac{1}{L}\Sigma_A^2)\hat{w}^\kappa+\frac{1}{L}\Sigma_A^2E\hat{x}^\kappa-\frac{1}{L}\Sigma_A^2\alpha + \frac{1}{L}q\Big)\\
		s^{\kappa+1} &= P_{\mathcal{S}} \Big( (I-\frac{1}{L}\Sigma_V^2)\hat{s}^\kappa+\frac{1}{L}\Sigma_V^2N\hat{x}^\kappa+\frac{1}{L}v\Big)\quad,
	\end{aligned}\label{eq:update_component}
\end{equation}
where $P_{\mathcal{Y}}$, $P_{\mathcal{W}}$, $P_{\mathcal{S}}$ are the projections onto the subsets of $\mathcal{Z}$ of the respective variables $y$, $w$, $s$. We refer the reader to Section~\ref{sec:FISTA_instantiation} for the definition of the projection operator and sets $\mathcal{Y}, \mathcal{W}, \mathcal{Z}, \mathcal{Z}$; and to Section~\ref{sec:Formulation_final} for the definition of the different matrices. Below, we consider the update for $x^{\kappa+1}$ separately, as it follows a different structure, while the remaining ones are similar.

\subsubsection{Update for the first component}

Let us first reformulate the expression for $x^{\kappa+1}$ in \eqref{eq:update_component} as
\begin{equation}
	\begin{aligned}
		x^{\kappa+1}   = x^{\kappa}  &+ \frac{1}{L} \Big[A^T \Sigma_N^2 ( \hat{y}^{\kappa}  -A \hat{x}^{\kappa} ) +  N^T\Sigma_V^2(\hat{s}^{\kappa} - N\hat{x}^\kappa)+   E^T\Sigma_A^2(\hat{w}^{\kappa}+\alpha - E\hat{x}^{\kappa}) \Big], \\ 
	\end{aligned}\label{eq:updatex_init}
\end{equation}where we have joined terms relative to the same edge matrices. We shall go over each of the three terms to obtain an update for $x_i^{\kappa+1}(\tau)$.

Consider the term $A^T \Sigma_N^2 ( \hat{y}^{\kappa}  -A \hat{x}^{\kappa} ) $. Only for derivation purposes, we define an auxiliary variable 
$$  \theta_{ij} = \frac{1}{\sigma_{ij}^2}(y_{ij}^\kappa(\tau) - (x_i^{\kappa}(\tau)- x_j^{\kappa}(\tau)) $$ and we define $\theta ( \tau)=\{\theta_{ij}(\tau)\}_{i\sim j}$ and $\theta = \{\theta(\tau)\}_{t-T_0\leq \tau \leq t}$. It follows that 
$$ \theta =  \Sigma_N^2 ( \hat{y}^{\kappa}  -A \hat{x}^{\kappa} ) $$ and, consequently, we consider the term $A^T\theta$. Multiplication on the left by $A^T$ equates to summing all terms in $\theta$ related to each node $i$ at time $\tau$. Therefore, we define $C_{(i\sim j,i)}$ as the element of $C$ on row referring to edge $i\sim j$ and column referring to node $i$ and obtain the update relative to $x_i(\tau)$ as
\begin{equation}\begin{aligned}
		\sum_{i\sim j \in \mathcal{E}_i} &\frac{1}{\sigma_{ij}^2}(C_{(i\sim j,i)}(\tau)y_{ij}^\kappa(\tau) + x_j^{\kappa}(\tau) - x_i^{\kappa}(\tau) )  = \\&- x_i^{\kappa}(\tau) \sum_{i\sim j \in \mathcal{E}_i} \frac{1}{\sigma_{ij}^2} + \sum_{i\sim j \in \mathcal{E}_i} \frac{1}{\sigma_{ij}^2}(C_{(i\sim j,i)}(\tau)y_{ij}^\kappa(\tau) + x_j^{\kappa}(\tau) ).  \\
	\end{aligned}\label{eq:aux_update_1}\end{equation}

Now, consider the term $E^T\Sigma_A^2(\hat{w}^{\kappa}+\alpha - E\hat{x}^{\kappa})$. Following the same reasoning we obtain the update relative to $x_i(\tau)$ as

\begin{equation}\begin{aligned}  \sum_{k \in \mathcal{A}_i}&\frac{1}{\varsigma_{ik}^2}\big(w_{ik}^\kappa(\tau)+\alpha_k(\tau) - x_i^\kappa(\tau) \big) = \\&- x_i^\kappa(\tau)  \sum_{k \in \mathcal{A}_i}\frac{1}{\varsigma_{ik}^2}+ \sum_{k \in \mathcal{A}_i}\frac{1}{\varsigma_{ik}^2}\big(w_{ik}^\kappa(\tau)+\alpha_k(\tau)\big) .\end{aligned}\label{eq:aux_update_2}\end{equation}

Finally, consider the term $N^T\Sigma_V^2(\hat{s}^{\kappa} - N\hat{x}^\kappa) $. Recall that  $N= C_{\textrm{vel}} \otimes I_p$ and each node has an edge with its position at $t-1$ and $t+1$, except for the first and last instants of the time window. Since the edges do not connect different nodes, we can define $N_i$ as the sub-matrix of $N$, whose rows correspond to node $i$. Further defining $\textbf{s}_i$ and $\textbf{x}_i$ as the concatenation of $s_i(\tau)$ and $x_i(\tau)$ over the time window, we get $ N_i^T\Sigma_V^2(\textbf{s}_i^{\kappa} - N_i\textbf{x}_i^\kappa) $, as a vector for the entire time window relative to node $i$. We then define 
\begin{equation}\big[N_i^T\Sigma_V^2(\textbf{s}_i^{\kappa} - N_i\textbf{x}_i^\kappa) \big]_{\tau} \label{eq:aux_update_3}\end{equation}
as the element relative to instant $\tau$.

Since we have obtained a separate update for $x_i(\tau)$ for each of the terms in \eqref{eq:updatex_init}, (and $x^\kappa$ is trivially distributed), we get the complete update for  $x^{\kappa+1}_i(\tau)$ from \eqref{eq:aux_update_1}, \eqref{eq:aux_update_2} and \eqref{eq:aux_update_2} as
\begin{equation*}
	\begin{aligned}
		x^{\kappa+1}_i(\tau) 
		=\frac{1}{L} \Bigg[  &x_i^\kappa(\tau) \Big(L - \sum_{i\sim j \in \mathcal{E}_i} \frac{1}{\sigma_{ij}^2} -\sum_{k \in \mathcal{A}_i}\frac{1}{\varsigma_{ik}^2}\Big)   +   \sum_{i\sim j \in \mathcal{E}_i} \frac{1}{\sigma_{ij}^2}\Big(C_{(i\sim j,i)}(\tau)y_{ij}^\kappa(\tau) \\& + x_j^{\kappa}(\tau) \Big)  + \big[N_i^T\Sigma_V^2(\textbf{s}_i^{\kappa} - N_i\textbf{x}_i^\kappa) \big]_{\tau}    + \sum_{k \in \mathcal{A}_i}\frac{1}{\varsigma_{ik}^2}\big(w_{ik}^\kappa(\tau)+\alpha_k(\tau)\big)  \Bigg].\\
	\end{aligned}
\end{equation*}

For convenience, we shall write this update as
$$ x_i^{\kappa+1}(\tau) = \frac{1}{L}\Big[ F_1(\tau)+F_2(\tau)+F_3(\tau)+F_4(\tau)  \Big], $$ where
\begin{equation*}
	F_1(\tau) =  \hat{x}_i^\kappa(\tau) \Big(L - \sum_{i\sim j \in \mathcal{E}_i} \frac{1}{\sigma_{ij}^2} -\sum_{k \in \mathcal{A}_i}\frac{1}{\varsigma_{ik}^2}\Big),
\end{equation*}
\begin{equation*}
	F_2(\tau) = \sum_{i\sim j \in \mathcal{E}_i} \frac{1}{\sigma_{ij}^2}\Big(\hat{x}_j^{\kappa}(\tau)+C_{(i\sim j,i)}(\tau)\hat{y}_{ij}^\kappa(\tau)\Big),
\end{equation*}
\begin{equation*}
	F_3(\tau) = \sum_{k \in \mathcal{A}_i}\frac{1}{\varsigma_{ik}^2}(\hat{w}_{ik}^\kappa(\tau)+\alpha_k(\tau))
\end{equation*}
and
\begin{equation*}
	F_4(\tau)=\big[N_i^T\Sigma_V^2(\hat{\textbf{s}}_i^{\kappa} - N_i\hat{\textbf{x}}_i^\kappa) \big]_{\tau}.
\end{equation*}

\subsubsection{Update for the remaining components}For the update of $y,w,z$, note that we can reformulate the equations in \eqref{eq:update_component} as 
\begin{equation*}
	\begin{aligned}
		y^{\kappa+1}&= P_{\mathcal{Y}} \Bigg( \frac{1}{L} \big[ u + \Sigma_N^2A\hat{x}^\kappa   + (IL-\Sigma_N^2)\hat{y}^\kappa   \big]    \Bigg)\\
		w^{\kappa+1} &= P_{\mathcal{W}} \Bigg( \frac{1}{L} \Big[ q - \Sigma_A^2\alpha +\Sigma_A^2E\hat{x}^\kappa+ (IL-\Sigma_A^2)\hat{w}^\kappa   \Big]    \Bigg)\\
		s^{\kappa+1} &= P_{\mathcal{S}} \Bigg( \frac{1}{L} \Big[ v +\Sigma_V^2N\hat{x}^\kappa+ (IL-\Sigma_V^2)\hat{s}^\kappa   \Big]    \Bigg).\\
	\end{aligned}
\end{equation*} 

It follows that the update for $y_{ij}^{\kappa+1}(\tau)$ is given as
\begin{equation}
	\begin{aligned}
		y_{ij}^{\kappa+1}(\tau) =  &P_{\mathcal{Y}_{ij}(\tau)} \Bigg(\frac{1}{L} \Big[ \tilde{u}_{ij}(\tau) +\frac{L\sigma_{ij}^2 - 1}{\sigma_{ij}^2} \hat{y}_{ij}^\kappa(\tau) + \frac{C_{(i\sim j,i)}(\tau) }{\sigma_{ij}^2} \Big(\hat{x}_i^\kappa(\tau)-\hat{x}_j^\kappa(\tau)\Big)  \Big]    \Bigg),\\
	\end{aligned}\label{eq:up1}
\end{equation} the update for $w_{ik}^{\kappa+1}(\tau)$ as
\begin{equation}
	\begin{aligned}
		w_{ik}^{\kappa+1}(\tau) = &P_{\mathcal{W}_{ik}(\tau)} \Bigg( \frac{1}{L} \Big[ \tilde{q}_{ik}(\tau) +\frac{ L\varsigma_{ik}^2 -1}{\varsigma_{ik}^2} \hat{w}_{ij}^\kappa(\tau)  +\frac{1 }{\varsigma_{ik}^2} \Big(\hat{x}_i^\kappa(\tau)-a_k(\tau)\Big)  \Big]    \Bigg)\\
	\end{aligned}\label{eq:up2}
\end{equation} and the update for $s_{i}^{\kappa+1}(\tau)$ as
\begin{equation}
	\begin{aligned}
		s_{i}^{\kappa+1}(\tau) = &P_{\mathcal{S}_{i}(\tau)} \Bigg( \frac{1}{L} \Big[ \tilde{v}_{i}(\tau) +\frac{ L\sigma_{i}^2 -1}{\sigma_{i}^2} \hat{s}_{i}^\kappa(\tau) + \frac{1 }{\sigma_{i}^2} \Big(\hat{x}_i^\kappa(\tau)-\hat{x}_i^\kappa(\tau-1)\Big)  \Big]    \Bigg).\\
	\end{aligned}\label{eq:up3}
\end{equation}
Note that updates in \eqref{eq:up1}, \eqref{eq:up2} and \eqref{eq:up3} only require information from node $i$ or its neighbours.

\vfill

\section{Algorithm with Parameter Estimation}
\label{appendix_B}

Algorithm~\ref{alg:ALGEstimation} is the extension of Algorithm~\ref{alg:ALG1} with parameter estimation. Bold symbols denote the concatenation of respective variables in the following way
\begin{equation*}
	\begin{split}
		\textbf{d}(t) &= \{d_{ij}(t) :i \sim j \in \mathcal{E}\}\\
		\textbf{u}(t) &= \{ u_{ij}(t):i \sim j \in \mathcal{E}\}\\
		\textbf{r}(t) &= \{r_{ik}(t):i \in \mathcal{V},k \in \mathcal{A}\}\\
		\textbf{q}(t) &= \{ q_{ik}(t):i \in \mathcal{V},k \in \mathcal{A}\}\\
		\boldsymbol{\beta}(t) &= \{\beta_i(t) \in \mathcal{V} \}\\
		\boldsymbol{\sigma}(t) &= \{\sigma_{ij}(t):i \sim j \in \mathcal{E}\}\\
		\boldsymbol{\varkappa}(t) &= \{\varkappa_{ij}(t):i \sim j \in \mathcal{E}\}\\
		\boldsymbol{\varsigma}(t) &= \{\sigma_{ik}(t):i \in \mathcal{V},k \in \mathcal{A}\}\\
		\boldsymbol{\lambda}(t) &= \{\lambda_{ij}(t): i \in \mathcal{V},k \in \mathcal{A}\}\\
		\boldsymbol{\sigma_v}(t) &= \{\sigma_{i}(t): i \in \mathcal{V}\}\\
		\boldsymbol{\varkappa_v}(t) &= \{\varkappa_{i}(t):i \in \mathcal{V}\}.
	\end{split}
\end{equation*}

\begin{algorithm}[H]
	\begin{algorithmic}[1]
		\REQUIRE{ $L$;\\ $z^0=(x^0,y^0,w^0,s^0)$,
			$\textbf{d}(t), \textbf{u}(t)$, $\textbf{r}(t), \textbf{q}(t) $, $\boldsymbol{\beta}(t)$, $\boldsymbol{\sigma}^0,\boldsymbol{\varkappa}^0 ,\boldsymbol{\varsigma}^0 ,\boldsymbol{\lambda}^0,\boldsymbol{\sigma_v}^0,\boldsymbol{\kappa_v}^0 $}; \\
		
		\ENSURE{ $\hat{x}(t)=\{\hat{x}_i(t)\}$}
		\STATE{$\sigma_{ij}' = \varsigma_{ik}' = \gamma_{ij}= \delta_{ik}= 0 $}
		\FORALL{t} 
		\STATE{Steps 2 to 26 of Algorithm 1}
		{\begin{equation*}\begin{split}
					\hat{\textbf{z}}(t) =\textrm{Algorithm\_1}\Big(& \hat{\textbf{z}}(t-1), \textbf{d}(t), \textbf{r}(t), \textbf{u}(t), \textbf{q}(t), \boldsymbol{\beta}(t), \boldsymbol{\sigma}(t), \boldsymbol{\varsigma}(t),\boldsymbol{\kappa}(t), \\& \boldsymbol{\lambda}(t),  \boldsymbol{\sigma_v}(t), \boldsymbol{\kappa_v}(t), L(t)\Big)
				\end{split}
		\end{equation*}}	
		\FORALL{node $i$}
		\FORALL{$i\sim j \in \mathcal{E}_i$}
		\STATE{$\hat{d}_{ij}(t) = \|\hat{x}_i(t)-\hat{x}_j(t)\|$, $\hat{u}_{ij}(t) = \frac{\hat{x}_i(t)-\hat{x}_j(t)}{\|\hat{x}_i(t)-\hat{x}_j(t)\|}$; }
		\STATE{$\sigma_{ij}' = \sigma_{ij}'+ (d_{ij}(t)-\hat{d}_{ij}(t))^2$, $\gamma_{ij} =\gamma_{ij} + \frac{u_{ij}(t)-\hat{u}_{ij}(t)}{\| u_{ij}(t)-\hat{u}_{ij}(t)\|} $;}
		\STATE{$\hat{\sigma}_{ij}^2(t+1) = \frac{1}{t}\sigma_{ij}'$, $\bar{\gamma}_{ij}(t) = \frac{\|\gamma_{ij}\|}{t}$}
		\STATE{$\hat{\varkappa}_{ij}(t+1) = \frac{\bar{\gamma}_{ij}(t)(p-\bar{\gamma}_{ij}^2(t))}{1-\bar{\gamma}_{ij}(t)^2}$}
		\ENDFOR
		\FORALL{$k \in \mathcal{A}_i$ }
		\STATE{$\hat{r}_{ik}(t) = \|\hat{x}_i(t)-a_k(t)\|$, $\hat{q}_{ik}(t) = \frac{\hat{x}_i(t)-a_k(t)}{\|\hat{x}_i(t)-a_k\|}$;}
		\STATE{$\varsigma_{ik}'=\varsigma_{ik}'+ (r_{ik}(t)-\hat{r}_{ik}(t))^2$, $\delta_{ik}=\delta_{ik}+\frac{q_{ik}(t)-\hat{q}_{ik}(t)}{\|q_{ik}(t)-\hat{q}_{ik}(t)\|}$;}
		\STATE{$\hat{\varsigma}_{ik}^2(t+1) = \frac{1}{t}\varsigma_{ik}'$, $\bar{\delta}_{ik}(t) = \frac{\| \delta_{ik}\|}{t}$}
		\STATE{$\hat{\lambda}_{ik}(t+1) = \frac{\bar{\delta}_{ij}(t)(p-\bar{\delta}_{ij}^2(t))}{1-\bar{\delta}_{ij}(t)^2}$}
		\ENDFOR
		\IF{$t>7$}
		\STATE{$\hat{v}_i(t)\hat{V}_i(t)\Delta T =$}
		\STATE{$=\frac{ 5(\hat{x}_i(t-3) - \hat{x}_i(t-5))}{32} + \frac{ 4(\hat{x}_i(t-2) - \hat{x}_i(t-6))}{32}+\frac{ (\hat{x}_i(t-1) - \hat{x}_i(t-7))}{32}$}
		\STATE{$V_i(t)v_i(t) = \beta(t) $}
		\STATE{$\sigma_{i}'= \sigma_{i}'+(V_i(t)-\hat{V}_i(t))^2$,  $\gamma_i= \gamma_i + (v_i(t)-\hat{v}_i(t))$;}
		\STATE{$\hat{\sigma}_{i}^2(t+1) = \frac{1}{t}\sigma_{i}'$, $\bar{\gamma}_{i}(t) = \frac{\|\gamma_i\|}{t}$}
		\STATE{$\hat{\varkappa}_{i}(t+1) = \frac{\bar{\gamma}_{i}(t)(p-\bar{\gamma}_{i}^2(t))}{1-\bar{\gamma}_{i}(t)^2}$}
		\ENDIF
		\ENDFOR
		\RETURN $\hat{\textbf{x}}(t)$
		\ENDFOR
	\end{algorithmic}
	\caption{Parameter Estimation}
	\label{alg:ALGEstimation}
\end{algorithm}

\section{Additional Simulation details}
\label{appendix:Exp}

\subsection{Details on simulation settings}

\paragraph{Trajectories}{ The lawnmower trajectory starts with all nodes separated by 15m on both axes and then travelling at a speed of \mbox{$1m/s$}. In the helix trajectory, the two nodes are separated by 10m and separated from the closest anchors by 8m and from the furthest by 15m. }

\begin{figure}[!htb]
	\centering
	\includegraphics[width =\linewidth]{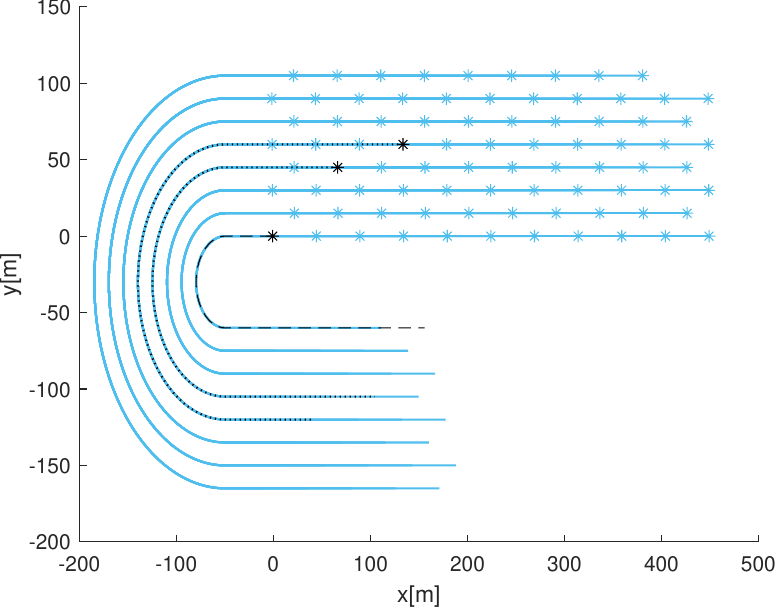}
	\caption{{Example of a trajectory with 80 nodes used for large scale experiments. Nodes are depicted in blue stars and the corresponding trajectories with blue solid lines. Anchors are depicted with black stars and their trajectories with dashed black lines. Nodes describe a partial lap trajectory.}}
	\label{fig:largescale}
\end{figure}

\paragraph{Large scale setting}{ For large scale experiments a partial lap trajectory is considered \mbox{(Figure~\ref{fig:largescale}}). The nodes are initially placed in the positions indicated with stars and then three random nodes are selected as anchors. A disk radius of $45$m is used to define neighbours for each of the nodes and the network topology is assumed constant throughout the trajectory. All nodes travel at a speed of $1$m/s.}

\subsection{Mean Positioning Errors}
Table~\ref{tab:errors} contains the MPE for the three trajectories and all the considered methods.

\begin{table}[!htb]
	\begin{tabular}{llllll}
		& \textbf{EKF}    & \textbf{Static} & \textbf{Ours Opt.} & \textbf{Ours Est.}  & \textbf{Ours No Est.}\\ ine
		\textbf{Lawnmower} & 0.1814 & 0.3675 & 0.1813                 & 0.1894          & 0.2113                 \\
		\textbf{Lap}       & 0.1792 & 0.3674 & 0.1817                 & 0.1903          & 0.2131                 \\
		\textbf{Helix}     & 0.1635 & 0.3070 & 0.1665                 & 0.1789           & 0.1973               
	\end{tabular}\caption{{MPE for each trajectory and all methods.}}\label{tab:errors}
\end{table}

\subsection{Different speeds}

\begin{figure}[!htb]
	\centering
	\includegraphics[width =0.5\linewidth]{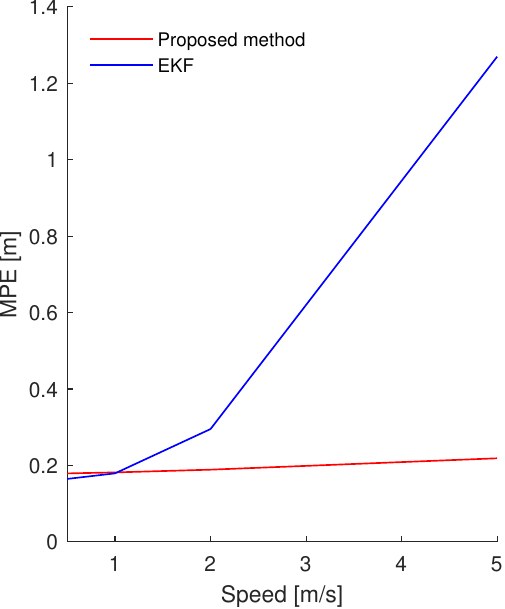}
	\caption[]{{MPE for different speed values during a lap trajectory.} }
	\label{fig:speeds}
\end{figure}

{We compare our method and EKF during a lap trajectory for different values of speed (\mbox{Figure	~\ref{fig:speeds}}). However, we note that in this case EKF was not tuned for the different speed measurements, except for the default value of $1$m/s. While the positioning error increases for both methods as the speed increases, our proposed solution suffers an increase of less than 0.1m, while EKF increases by one order of magnitude.}
	
\end{document}